\documentclass{article}%
\usepackage{amsmath}
\usepackage{amsfonts}
\usepackage{amssymb}
\usepackage{graphicx}
\usepackage{color}%
\newtheorem{theorem}{Theorem}

\newtheorem{definition}[theorem]{Definition}

\newtheorem{lemma}[theorem]{Lemma}

\newtheorem{proposition}[theorem]{Proposition}

\newenvironment{proof}[1][Proof]{\noindent\textbf{#1.} }{\ \rule{0.5em}{0.5em}}

\begin{document}

\title{Robust location estimation with missing data }
\author{Mariela Sued (msued@dm.uba.ar) and V\'{\i}ctor J. Yohai (vyohai@uolsinectis.com.ar)\\University of Buenos Aires and CONICET}
\date{~}
\maketitle
\begin{abstract}
In a missing-data setting, we have a sample in which a vector of
explanatory variables $\mathbf{x}_{i}$ is observed for every subject
$i$, while scalar outcomes $y_{i}$ are missing by happenstance on
some individuals. In this work we propose robust estimates of the
distribution of the responses assuming missing at random (MAR) data,
under a semiparametric regression model. Our approach allows the
consistent estimation of any weakly continuous functional of the
response's distribution. In particular, strongly consistent
estimates of any continuous location functional, such as the median
or MM functionals, are proposed. A robust fit for the regression
model combined with the robust properties of the location functional
gives rise to a robust recipe for estimating the location parameter.
Robustness is quantified through the breakdown point of the proposed
procedure. The asymptotic distribution of the location estimates is
also derived.

\end{abstract}

\section{Introduction}

Suppose we have a sample of a population, such that for every subject $i$ in
the sample we observe a vector of explanatory variables $\mathbf{x}_{i}$ while
a scalar response $y_{i}$ is missing by happenstance on some individuals.\ A
classical problem is to construct consistent estimators for the mean value of
the response based on the observed data. In order to identify the parameter of
interest in terms of the distribution of observed data, missing at random
(MAR) is assumed.

This hypothesis establishes that the value of the response does not
provide additional information, on top of that given by the
explanatory variables, to predict whether an individual will present
a missing response (see Rubin \cite{Rubin}). To be more rigorous,
let us introduce a binary variable $a_{i}$ such that $a_{i}=1$
whenever the response is observed for subject $i$. In this way,
MAR states that%
\begin{equation}
\mathrm{P}(a_{i}=1|\mathbf{x}_{i},y_{i})=\mathrm{P}(a_{i}=1|\mathbf{x}_{i}).
\label{MAR}%
\end{equation}
Under this condition, if $\mathrm{P}(a_{i}=1|\mathbf{x}_{i})>0$, we have that%
\begin{equation}
\mathrm{E}[y_{i}]=\mathrm{E}\left[  \frac{a_{i}y_{i}}{\pi(\mathbf{x}_{i}%
)}\right]  , \label{IPW}%
\end{equation}
where $\pi(\mathbf{x}_{i})=$\textrm{$P$}$(a_{i}=1|\mathbf{x}_{i})$, and
identifiability of $\mathrm{E}[y_{i}]$ holds. One approach to estimate
consistently $\mathrm{E}[y_{i}]$, called inverse probability weight (IPW), is
based on (\ref{IPW}) and requires to estimate the propensity score function
$\pi(\mathbf{x)}$. Then, the estimate of $\mathrm{E}[y_{i}]$ can be obtained
replacing in (\ref{IPW}) $\pi(\mathbf{x}_{i})$ by its estimate and the
expectation by its empirical version. MAR also implies that the conditional
distribution of the responses given the vector of explanatory variables
remains the same, regardless of the fact that the response is also observed:
$y_{i}|\mathbf{x}_{i}\sim y_{i}|\mathbf{x}_{i},a_{i}=1$. Then $\mathrm{E}%
[y_{i}|\mathbf{x}_{i}]=\mathrm{E}[y_{i}|\mathbf{x}_{i},a_{i}=1]$. Since
$\mathrm{E}[y_{i}]=\mathrm{E}\left[  \mathrm{E}[y_{i}|\mathbf{x}_{i}]\right]
$, a second approach to estimate $\mathrm{E}[y_{i}]$ is based \ on a
regression model (parametric or nonparametric ) for $\mathrm{E}[y_{i}%
|\mathbf{x}_{i}]=g(\mathbf{x}_{i})$, which is fitted using only the
individuals for whom the response is observed. \ Then a second
estimate for $\mathrm{E}[y_{i}]$ is obtained by averaging
$\widehat{g}(\mathbf{x}_{i})$ over the whole sample, where
$\widehat{g}$ is an estimate of $g$. There is third approach (doubly
protected) that \ postulates models for $\pi (\mathbf{x)}$ and
$g(\mathbf{x)}$ \ and obtains a consistent estimate of
$\mathrm{E}[y_{i}]$ if at least one of the two models is  correct. A
recent survey and discussion on these three approaches can be found
in Kan and Schafer \cite{KanSchafer}\ and Robins, Sued, Lei-Gomez
and Rotnitzky \cite{Robins}.

As it is well known, the mean is not a robust location parameter,
i.e., a small change in the population distribution may have a large
effect on this parameter. As a consequence of this, the mean does
not admit consistent non-parametric robust estimates, except when
strong properties on the distribution are assumed, as for example
symmetry. For this reason, to introduce robustness in the present
setting, we start by reformulating the statistical object of
interest: instead of estimating the mean value of the response, we
look for consistent estimates of \ $T_{L}(F_{0}),$ where $T_{L}$ is
a robust location functional and $F_{0}$ is the distribution of
$y_{i}$. Bianco, Boente, Gonzalez-Manteiga and Perez-Gonzalez
\cite{Bianco} used this approach to obtain robust and consistent
estimates of an M location parameter of the distribution of $y_{i}$.
In their treatment they assumed a partially linear model to describe
the relationship between $y_{i}$ and $\mathbf{x}_{i}$, and also that
the distributions of the \ response $y_{i}$ and of the regression
error under the true model are both symmetric.

In this paper we introduce a new estimate of any continuous location
functional assuming that the relation between $y_{i}$ and
$\mathbf{x}_{i}$ is given by means of a semiparametric regression
model. We show that once the regression model is fitted using robust
estimates, we can define a consistent estimate of the distribution
function of the response. Then, any parameter of the response
distribution defined throughout a weak continuous functional, may be
also consistently estimated by evaluating the functional at the
estimated distribution function. The consistency of this procedure
does not require the symmetry assumptions used by Bianco et al.
\cite{Bianco}.

A robust fit for the regression model combined with the robust properties of
the location functional to be considered, gives rise to a robust recipe for
estimating the location parameter. Robustness is quantified looking at
breakdown point\ of the proposed procedure. In particular our results can be
applied when the location functional is the median or an MM location functional.

The proposed procedure may be considered as a robust extension of the second
approach described above for estimating $\mathrm{E}[y_{i}]$. We have not found
a way to robustify the approaches that use the propensity score $\pi
(\mathbf{x)}$. The main difficulty in such cases is to obtain a consistent
procedure avoiding the assignment of very large weights to those observations
with very small $\pi(\mathbf{x}_{i}\mathbf{)}$.

This work is organized as follows. In Section 2 we formalize the problem of
the robust estimation of a location parameter with missing data. We propose a
family of procedures which depend on the location functional to be estimated
and also on the robust regression estimate for the parameter of the regression
model postulated to describe the relationship between $\mathbf{x}_{i}$ and
$y_{i}$. In Section 3 we show that, under some assumptions on the location
functional and the regression estimate, the proposed estimates are strongly
consistent and asymptotically normal. In Section 4 we study the breakdown
point of the proposed estimates. In Section 5 we show that when the location
and regression estimates are of MM type, the assumptions that guarantee
consistency and asymptotic normality of the proposed estimates are satisfied.
In Section 6 we present the results of a Monte Carlo study which shows that
the proposed estimates are highly efficient under Gaussian errors and highly
robust under outlier contamination. Proofs are presented in the Appendix.

\section{Notation and Preliminaries}

We first introduce some notation. Henceforth $\mathrm{E}_{G}[h(\mathbf{z})]$
and $\mathrm{P}_{G}\left(  A\right)  $ will respectively denote the
expectation of $h(\mathbf{z})$ and the probability that $\mathbf{z}\in A,$
when $\mathbf{z}$ is distributed according to $G$. If $\mathbf{z}$ has
distribution $G$ we write $\mathbf{z}\sim G$ or $\mathcal{D}\left(
\mathbf{z}\right)  =G.$ Weak convergence of distributions, convergence in
probability and convergence in distribution of random variables or vectors are
denoted by $G_{n}\rightarrow_{w}G,$ $\mathbf{z}_{n}\rightarrow_{p}\mathbf{z}$
and $\mathbf{z}_{n}\rightarrow_{d}\mathbf{z},$ respectively. \ By an abuse of
notation, we will write $\mathbf{z}_{n}\rightarrow_{d}G$ to denote
$\mathcal{D}\left(  \mathbf{z}_{n}\right)  \rightarrow_{w}G.$ We use
$o_{P}(1)$ to denote any sequence that converges to zero in probability. The
complement and the indicator of the set $A$ are denoted by $A^{c}$ and
$\mathbf{1}_{A},$ respectively. The scalar product of vectors $\mathbf{a,b}%
\in{\mathbb{R}}^{s}$ is denoted by $\mathbf{a}^{\prime}\mathbf{b}$.
$\mathbb{R}_{+}$ denotes the set of positive real numbers.

Along this paper \ we use the expression \textit{empirical distribution \ } of
$\ $a$\ $sequence on $\ n$ points $\mathbf{z}_{1},\mathbf{z}_{2}%
,...,\mathbf{z}_{n}$ in $\mathbb{R}^{k}$\textit{\ } to denote the function
$F_{n}:\mathbb{R}^{k}\rightarrow\lbrack0,1]$ such that given $\mathbf{z}%
\in\mathbb{R}^{k},$ $F_{n}(\mathbf{z})$ $=m/n,$ where $m$ is the number of
points $\mathbf{z}_{i}$ such that \ all \ its coordinates $\ $ are smaller or
equal than the corresponding ones of $\mathbf{z}.$

\subsection{Describing our setting: the data, the problem and the model}

Throughout this work, \ we have \ a random sample of $\ n$ subjects \ and for
each subject $i$ \ in the sample, $1\leq i\leq n$, a vector of explanatory
variables $\mathbf{x}_{i}$ is always observed, while the response $y_{i}$ is
missing on some subjects. Let $a_{i}$ be the indicator of whether $y_{i}$ is
observed at subject $i$: $a_{i}=1$ if $y_{i}$ is observed and $a_{i}=0$ if it
is not.

We will be concerned with the estimation of a location \ functional at the
distribution of the response. A location functional $T_{L}$, defined on a
class of univariate distribution functions $\mathcal{G}$, assigns to each
$F\in$ $\mathcal{G}$ a real number $\ T_{L}(F)$ satisfying $T_{L}%
(F_{ay+b})=aT_{L}(F_{y})+b,$ where $F_{y}$ denotes the distribution of the
random variable $y$.

Example of locations functionals are the mean and median. Another
important class of location functionals that includes the mean and
median and other robust estimates \ is the class of M location
functionals. \ This class also includes S and MM estimators that
will be described in Section 4. We should also mention the class of
L location functionals, see e.g. Chapter 2 of Maronna, Martin and
Yohai \cite{MaronnaMartin and Yohai}, but we do not study this class
in this work.

A functional $T$ is said to be weakly continuous at $F$ if \ given a sequence
$\{F_{n}\}$ of distribution functions that converges weakly to $F$
$(F_{n}\rightarrow_{w}F$), then $T(F_{n})\rightarrow T(F)$. In order to obtain
a consistent estimate of a location parameter defined by means of a weakly
continuous functional, \ it is sufficient to \ have \ a sequence \ of
estimates $\widehat{F}_{n}$ \ such that \ converges weakly to the distribution
of the $y_{i}$'s.

To be more precise, denote by $F_{0}$ the distribution of the outcomes $y_{i}
$. \ Let $\ T_{L}$ be a weakly continuous location functional at $F_{0}$ . We
are interested in estimating
\[
\mu_{0}=T_{L}(F_{0}).
\]

We assume a semiparametric regression model
\begin{equation}
y_{i}=g(\mathbf{x}_{i},\mathbf{\beta}_{0})+u_{i},1\leq i\leq n,
\label{modelogral}%
\end{equation}
with $y_{i},u_{i}\in{\mathbb{R}}$, $\mathbf{x}_{i}\in{\mathbb{R}}^{p}$,
$u_{i}$ independent of $\mathbf{x}_{i}$, \textbf{$\beta$}$_{0}\in
B\subset{\mathbb{R}}^{q}$, $g:{\mathbb{R}}^{p}\times B\rightarrow{\mathbb{R}}%
$. Furthermore, in order to guarantee the MAR condition, we assume that
$u_{i}$ is independent of $(\mathbf{x}_{i},a_{i})$. We denote by $Q_{0}$ and
$\mathrm{K}_{0}$ the distributions of $\mathbf{x}_{i}$ and $u_{i}$, respectively.

To identify \textbf{$\beta$}$_{0}$, \ without assuming that either (i)
$\mathrm{K}_{0}$ is symmetric around $0$ or (ii) \ $\mathrm{K}_{0}$ satisfies
a centering condition, (as, e.g. , E$_{K_{0}}u=0)$ we assume that
\begin{equation}
\mathrm{P}_{Q_{0}}\left(  g(\mathbf{x},\mathbf{\beta}_{0})=g(\mathbf{x}%
,\mathbf{\beta})+\alpha\right)  <1 \label{IDCOND}%
\end{equation}
for all \textbf{$\beta$}$\neq$\textbf{$\beta$}$_{0}$, for all $\alpha$. \ This
condition requires that in case there is an intercept, it will be included in
the error term $u_{i}$ instead as of a parameter of the regression function
$g(\mathbf{x},$\textbf{$\beta$}$)$. For linear regression we have
$g(\mathbf{x},$\textbf{$\beta$}$)=$\textbf{$\beta$}$^{\prime}\mathbf{x}$ and
then this condition means that the vector $\mathbf{x}_{i}$ \ \ is not
concentrated on any hyperplane.

\subsection{The proposal}

Recall that $\mathrm{K}_{0}$ denotes de distribution of $u_{i}$ and let
$R_{0}$ denotes the distribution of $g(\mathbf{x}_{i},\beta_{0})$.
Independence between $\mathbf{x}_{i}$ and $u_{i}$ guarantees that $F_{0}$ is
the convolution between $R_{0}$ and $K_{0}$. Then by convoluting consistent
estimators $\widehat{R}_{n}$ and $\widehat{K}_{n}$ of each of these
distributions, we get a consistent estimator for $F_{0}$.

In order to estimate $R_{0}$ and $K_{0}$ we need to have a robust \
and strongly consistent estimator $\widehat{\mathbf{\beta}}_{n}$ of
\textbf{$\beta $}$_{0}$. This estimator may be, for example, an \ S
estimate (see Rousseeuw and Yohai \cite{RousseeuwYohai}) or an
MM-estimate (see Yohai \cite {Yohai87}). Since $u_{i}$ is
independent of $a_{i}$, $\widehat{\mathbf{\beta}}_{n}$ may be
obtained by a robust fit of the model using the data for which $\
y_{i}$ is observed: \ i.e., using the observations
$(\mathbf{x}_{i},y_{i})$ with $a_{i}=1$. Let
$\widehat{R}_{n}$ be the empirical distribution of $g(\mathbf{x}_{j}%
,\widehat{\mathbf{\beta}}_{n})$, $1\leq j\leq n$ defined by%
\begin{equation}
\widehat{R}_{n}=\frac{1}{n}\,\sum_{j=1}^{n}\delta_{g(\mathbf{x}_{j}%
,\widehat{\mathbf{\beta}}_{n})}, \label{Rn}%
\end{equation}
where $\delta_{s}$ denotes the point mass distribution \ at $s$. \

Let $A=\{i:a_{i}=1\}$ and $m=\#A$. For $\ i\in A$ consider
\[
\widehat{u}_{i}=y_{i}-g(\mathbf{x}_{i},\widehat{\mathbf{\beta}}_{n}).
\]
The estimator $\widehat{K}_{n}$ of $K_{0}$ is defined as the empirical
distribution of $\{\widehat{u}_{i}:i\in A\}$:
\begin{equation}
\widehat{K}_{n}=\frac{1}{m}\,\sum_{i\in A}\delta_{\widehat{u}_{i}}=\frac
{1}{\sum_{i=1}^{n}a_{i}}\,\sum_{i=1}^{n}a_{i}\delta_{\widehat{u}_{i}}.
\label{Fn}%
\end{equation}

Then, we estimate $F_{0}$ by $\widehat{F}_{n}=\widehat{R}_{n}\ast\widehat
{K}_{n}$, where $\ast$ denotes convolution. Note that $\widehat{R}_{n}%
\ast\widehat{K}_{n},$ is the empirical distribution of the $nm$ points%
\[
\widehat{y}_{ij}=g(\mathbf{x}_{j},\widehat{\mathbf{\beta}}_{n})+\widehat
{u}_{i},\text{ }1\leq j\leq n,\text{ }i\in A,
\]
and therefore we can also express $\widehat{F}_{n}$ as%
\begin{equation}
\widehat{F}_{n}=\frac{1}{nm}%
{\displaystyle\sum\limits_{i\in A}}
\sum_{j=1}^{n}\delta_{\widehat{y}_{ij}}=\frac{1}{n\,\sum_{i=1}^{n}a_{i}}%
{\displaystyle\sum\limits_{i\in A}}
\sum_{j=1}^{n}\delta_{\widehat{y}_{ij}}. \label{Gnsom}
\end{equation}
Finally, we estimate $\mu_{0}$ by
\begin{equation}
\widehat{\mu}_{n}=T_{L}(\widehat{F}_{n}). \label{OUREST}%
\end{equation}

Since we have assumed weak continuity of $T_{L}$ at $F_{0}$, in order to prove
that $\widehat{\mu}_{n}$ is \ a strongly consistent estimate of $\mu_{0}$ we
only need to prove that $\widehat{F}_{n}\rightarrow_{w}F_{0}$ a.s. Observe
that
\[
\mathrm{E}_{\widehat{F}_{n}}h(y)=\frac{1}{nm}%
{\displaystyle\sum\limits_{i\in A}}
\sum_{j=1}^{n}h(\widehat{y}_{ij}).
\]
The right hand \ side of this equation \ was proposed by M\"{u}ller
\cite{Muller} to estimate $\mathrm{E}_{F_{0}}h(y).$

\section{Consistency and asymptotic distribution\label{results}}

Let $(\mathbf{x}_{i},y_{i})$ and $\ u_{i}$ satisfy model (\ref{modelogral}),
with $u_{i}$ independent of $(\mathbf{x}_{i},a_{i})$. Denote by $\mathrm{G}%
_{0}$, $\mathrm{Q}_{0}$ and $\mathrm{K}_{0}$ the distributions of
$(\mathbf{x}_{i},y_{i})$, $\mathbf{x}_{i}$ and $u_{i}$, respectively, and
denote by $\mathrm{G}_{0}^{\ast}$ and $\mathrm{Q}_{0}^{\ast}$ the
distributions of $(\mathbf{x}_{i},y_{i})$ and $\mathbf{x}_{i}$ conditioned on
$a_{i}=1$, respectively.

The MAR condition implies that under $G_{0}^{\ast}$ model (\ref{modelogral})
is still satisfied with $\ \mathbf{x}_{i}^{\ast}$ and $u_{i}^{\ast}$
independent, $\mathbf{x}_{i}^{\ast}$ with distribution $Q_{0}^{\ast}$ and
$\ u_{i}^{\ast}$ with distribution $K_{0}.$ We also assume\ that the
regression function $\ g$ satisfies following assumption:

\textbf{A0 } $g(\mathbf{x},$\textbf{$\beta$}$)$ is twice continuously
differentiable with respect to \textbf{$\beta$} and there exists $\delta>0$
such that%
\begin{equation}
\mathrm{E}_{Q_{0}}\sup_{\left\Vert \mathbf{\beta}-\mathbf{\beta}%
_{0}\right\Vert \leq\delta}\left\Vert {\dot{g}}(\mathbf{x}_{1},\mathbf{\beta
})\right\Vert ^{2}<\infty\text{ }\mathrm{and}\ \ \mathrm{E}\,_{Q_{0}}%
\sup_{\left\Vert \mathbf{\beta}-\mathbf{\beta}_{0}\right\Vert \leq\delta
}\left\Vert {\ddot{g}}(\mathbf{x}_{1},\mathbf{\beta})\right\Vert <\infty,
\label{A0}%
\end{equation}

where $\dot{g}(\mathbf{x},$\textbf{$\beta$}$)$ and $\ddot{g}\left(
\mathbf{x},\mathbf{\beta}\right)  $ denote the vector of first derivatives and
the matrix of second derivatives of $g$ respect to \textbf{$\beta$}, respectively.

In order to prove the consistency and the asymptotic normality of
$\widehat{\mu}_{n}$ the following assumptions on $\widehat{\mathbf{\beta}}%
_{n}$ and $T_{L}$ are required.

\textbf{A1} $\{\widehat{\mathbf{\beta}}_{n}\}$ is strongly consistent for
\textbf{$\beta$}$_{0}$.

\textbf{A2 }The regression estimate $\widehat{\mathbf{\beta}}_{n}$ satisfies
\begin{equation}
\sqrt{n}(\widehat{\mathbf{\beta}}_{n}-\mathbf{\beta}_{0})=\frac{1}{n^{1/2}%
}\sum_{i=1}^{n}a_{i}I_{R}(\mathbf{x}_{i},y_{i})+o_{P}(1), \label{TRexpansion}%
\end{equation}
for some function $I_{R}(\mathbf{x},u)$ with $\mathrm{E}a_{i}I_{R}%
(\mathbf{x_{i}},y_{i})=0$ and finite second moments.

\textbf{A3} $T_{L}$ is weakly continuous at $F_{0}$.

\textbf{A4 } The following expansion holds:
\begin{equation}
\sqrt{n}\left(  T_{L}(\widehat{F}_{n})-T_{L}(F_{0})\right)  =\sqrt
{n}\mathrm{E}_{\widehat{F}_{n}}I_{L}(y)\,+\,o_{P}(1), \label{TLexpansion}%
\end{equation}
for some differentiable \ function $I_{L}(y)$ with $\mathrm{E}_{F_{0}}%
I_{L}(y)=0$, $\mathrm{E}_{F_{0}}I_{L}^{2}(y)<\infty$ and $|I_{L}^{\prime}(y)|$ bounded.

It can be shown that when expansion (\ref{TLexpansion}) holds, $I_{L}$ is
given by the influence function (as defined by Hampel (1974)) of $T_{L}$ at
$F_{0}$. When $\widehat{\mathbf{\beta}}_{n}$ is obtained using a regression
functional, a similar statement holds.

The following Theorem shows the consistency of $\widehat{\mu}_{n}%
=T(\widehat{F}_{n}).$

\begin{theorem}
\label{lemaconvdebil1} Let $\widehat{F}_{n}$ be defined as in
(\ref{Gnsom}) and assume that A1 holds. Then (a)
$\{\widehat{F}_{n}\}$ converges weakly to $F_{0}$ a.s., i.e.,
\[
\mathrm{P}(\widehat{F}_{n}\rightarrow_{w}F_{0})=1.
\]
(b) Assume\ also that A3 holds; then $\widehat{\mu}_{n}=T_{L}(\widehat{F}%
_{n})$ converges a.s. to $\mu_{0}=T_{L}(F_{0})$.
\end{theorem}

In order to find \ the asymptotic distribution of $\widehat{\mathbf{\mu}}%
_{n}\ $, consider
\begin{align*}
\eta &  =\mathrm{E}a_{i},\mathbf{c}=\mathrm{E}\left[  a_{1}\,I_{L}^{\prime
}(y_{1}-g({\mathbf{\beta}_{0}},\mathbf{x}_{1})+{g}({\mathbf{\beta}_{0}%
},\mathbf{x}_{2}))\,\{\dot{g}({\mathbf{\beta}_{0}},\mathbf{x}_{2})-\dot
{g}({\mathbf{\beta}_{0}},\mathbf{x}_{1})\}\right]  ,\\
e(\mathbf{x}_{i},u_{i},a_{i})  &  =\mathrm{E}\left[  a_{i}I_{_{T_{L},F_{0}}%
}(u_{i}+g(\mathbf{x}_{j},\mathbf{\beta}_{0}))|u_{i},a_{i}\right]
=a_{i}\mathrm{E}\left[  I_{_{T_{L},F_{0}}}(u_{i}+g(\mathbf{x}_{j}%
,\mathbf{\beta}_{0}))|u_{i},a_{i}\right]  ,\\
f(\mathbf{x}_{j}\mathbf{)}  &  \mathbf{=}\mathrm{E}\left[  a_{i}%
I_{_{T_{L},F_{0}}}(u_{i}+g(\mathbf{x}_{j},\mathbf{\beta}_{0}))|\mathbf{x}%
_{j}\right]  ,\\
\tau^{2}  &  =\frac{1}{\eta^{2}}\mathrm{E}\left[  \left\{  e(\mathbf{x}%
_{i},u_{i},a_{i})+f(\mathbf{x}_{i})+a_{i}{\mathbf{c}^{\prime}I_{R}%
(\mathbf{x}_{i},u_{i})}\right\}  ^{2}\right]  .
\end{align*}
Then, the following Theorem gives the asymptotic normality of the estimate
$\widehat{\mathbf{\mu}}_{n}$, defined in (\ref{OUREST}).

\begin{theorem}
\label{asym} Assume A0-A4. Then%
\begin{equation}
n^{1/2}(\widehat{\mu}_{n}-\mu_{0})\rightarrow_{d}N(0,\tau^{2}).
\label{Vnasnor}%
\end{equation}

\end{theorem}

\subsection{The median as location parameter}

The median is one of the most popular robust location parameters.
However, \ since this estimate does not satisfy A4, we cannot prove
its asymptotic normality using Theorem \ref{lemaconvdebil1}.\ In
this section, we will prove consistency and asymptotic distribution
for the median of $\widehat{F}_{n}$, defined at (\ref{Gnsom}),
assuming that A0 holds and that $\{\widehat{\beta }_{n}\}$ satisfy
A1 and A2.

The functional $T_{\text{med}}$ is defined by
\begin{equation}
T_{\text{med}}(F)=\arg\min_{\mu}\mathrm{E}_{F}|y-\mu|. \label{meddef}%
\end{equation}
When there are more than one value attaining the minimum, the
functional is defined \ by choosing \ any of them. We have the
following result, whose proof needs an extra argument to compensate
the absence of differentiability of $I_{T_{\text{med}},F_{0}}(y)$. \

\begin{theorem}
\label{median} \ Assume that $\mu_{0}=T_{\text{med}}(F_{0})$ is well
defined and let $\widehat{\mu}_{n}=T_{\text{med}}(\widehat{F}_{n})$.
 Suppose that $F_{0}$ is continuous and strictly
increasing at $\mu_{0}$. Then, \ (a) under A1 we have
$\widehat{\mu}_{n}\rightarrow\mu_{0}$ a.s.\newline(b )Assume A0-A2.
\ Assume also that $F_{0}$ and $K_{0}$ have continuous \ and
bounded densities $f_{0}$ $\ $and $k_{0}$ respectively, and that $f_{0}%
(\mu_{0})>0$. \ Then%
\begin{equation}
n^{1/2}(\widehat{\mu}_{n}-\mu_{0})\rightarrow_{d}N(0,\tau^{2}),
\label{Vnasnor1}%
\end{equation}
where $\tau^{2}$ is as in Theorem \ref{asym} with $\mathbf{c}$ replaced by%
\begin{equation}
\mathbf{c}^{\ast}=\frac{1}{\eta f_{0}(\mu_{0})}\mathrm{E}[a_{1}k_{0}%
(-g(\mathbf{x}_{2},\mathbf{\beta}_{0})+\mu_{0})\{\dot{g}(\mathbf{x}%
_{2},\mathbf{\beta}_{0})-\dot{g}(\mathbf{x}_{1},\mathbf{\beta}_{0}%
)\}]\ \label{cstar}%
\end{equation}
and $I_{T_{L},F_{0}}(y)$ replaced by
\[
I_{T_{\text{med}},F_{0}}(y)=\frac{\text{sign}(y-\mu_{0})}{2f_{0}(\mu_{0})}.
\]

\end{theorem}

\section{Breakdown point \label{secBDP}}

Consider first a dataset of $n$ complete observations $\mathbf{Z=\{z}%
_{1}\mathbf{,..,z}_{n}\mathbf{\}},$ where
$\mathbf{z}_{i}\in\mathbb{R}^{j}$, and let
$\widehat{\theta}_{n}(\mathbf{Z})$ be an estimate of a parameter
$\mathbf{\theta}\in\mathbb{R}^{k}$ defined on all possible datasets.
Donoho and Huber \cite{Dh} define the finite sample breakdown point
(FSBP) of $\widehat{\mathbf{\theta}}_{n}$ at $\ \mathbf{Z}$ by
\[
\varepsilon^{\ast}(\widehat{\mathbf{\theta}}_{n},\mathbf{Z})=\min\left\{
\frac{s}{n}:\sup_{\mathbf{Z}^{\ast}\in\mathcal{Z}_{s}}\Vert\widehat{\theta
}_{n}(\mathbf{Z}^{\ast})\Vert=\infty\right\}  ,
\]
where%
\[
\mathcal{Z}_{s}=\{\mathbf{Z}^{\ast}=\{\mathbf{z}_{1}^{\ast},...,\mathbf{z}%
_{n}^{\ast}\}:%
{\displaystyle\sum\limits_{i=1}^{n}}
I\{\mathbf{z}_{t}^{\ast}\neq\mathbf{z}_{i}\}\leq s\}.
\]
\ Then $\ \varepsilon^{\ast}$ is the minimum fraction of outliers that is
required to take the estimate\ beyond\ any bound.

Now, we extend \ the notion of FSBP to the present setting, where there are
missing data, as follows. \ Let\
\begin{equation}
\mathbf{W}=\{(\mathbf{x}_{1},y_{1},a_{1}),....(\mathbf{x}_{n},y_{n},a_{n})\}
\label{Wcom}%
\end{equation}
be the set of \ all observations\ and missingness\ \ indicators, and \ let
$A=\{i:1\leq i\leq n,\,a_{i}=1\}$, $m=\#A.$ Denote by $\mathcal{W}_{st}$ the
set of all samples \ obtained from $\mathbf{W}$ by replacing at most $t$
points by \ outliers, with at most $s$ of these replacement corresponding to
the non missing observations.\ Then $\mathbf{W}^{\ast}=\{(\mathbf{x}_{1}%
^{\ast},y_{1}^{\ast},a_{1}),....(\mathbf{x}_{n}^{\ast},y_{n}^{\ast},a_{n})\}$
belongs to $\mathcal{W}_{t,s}$ if
\[
\sum_{i\in A}I\{(\mathbf{x}_{i}^{\ast},y_{i}^{\ast})\neq(\mathbf{x}_{i}%
,y_{i})\}\ +\sum_{i\in A^{C}}I\{\mathbf{x}_{i}^{\ast}\neq\mathbf{x}_{i}\}\leq
t
\]
and
\[
\sum_{i\in A}I\{(\mathbf{x}_{i}^{\ast},y_{i}^{\ast})\neq(\mathbf{x}_{i}%
,y_{i})\}\ \leq s.
\]
Given an estimate $\widehat{\mu}_{n}$ of $\mu_{0},$ we define
\[
M_{ts}=\sup_{\mathbf{W}^{\ast}\in\mathcal{W}_{t,s}}\left\vert \widehat{\mu
}_{n}(\mathbf{W}^{\ast})\right\vert
\]
and
\[
\kappa(t,s)=\max\left(  \frac{t}{n},\frac{s}{m}\right)  .
\]
Then, we define the finite sample breakdown point \ (FSBP) of an estimate
\ $\widehat{\mu}_{n}$ at $\mathbf{W}$
\[
\varepsilon^{\ast}=\min\{\kappa(t,s):M_{ts}=\infty\}.
\]
Then $\ \varepsilon^{\ast}$ is the minimum fraction of outliers in the
complete sample or in the set of non missing \ observations that is required
to take the estimate \ beyond\ any bound.

In order to get a lower bound for the FSBP of the location estimate
$\widehat{\mu}_{n}$ introduced in (\ref{OUREST}), \ we need to define the
\emph{\ uniform asymptotic breakdown point} $\varepsilon_{U}^{\ast}$\ of
$T_{L}$ as follows:

\begin{definition}
Given a functional $T_{L}$, its \emph{uniform asymptotic breakdown point}
$\ $(UABP) $\varepsilon_{U}^{\ast}( T_{L})$ is defined as the supremum of all
$\varepsilon>0$ satisfying the following property: for all $M>0$ there exists
$K>0$ depending on $M$ so that
\begin{equation}
\mathrm{P}_{F}(|y|\leq M)>1-\varepsilon\text{ implies }|T_{L}(F)|<K.
\label{epsilonun}%
\end{equation}

\end{definition}

For any location functional $T_{L}$ we have that $\varepsilon_{U}^{\ast}%
(T_{L})\leq0.5.$ This is an immediate consequence of the following two facts:
(a) $\varepsilon_{A}^{\ast}$\ $($ $T_{L},F)$ $\leq0.5$, for all location
functionals $T_{L}$ and all $F$, where $\varepsilon_{A}^{\ast}$\ $($
$T_{L},F)$ is the asymptotic breakdown point of $T_{L}$ at the distribution
$F$, while (b) $\varepsilon_{U}^{\ast}(T_{L})$ $\leq\varepsilon_{A}^{\ast}%
$\ $($ $T_{L},F)$ for all $F$. In the case that $T_{L}$ is the median it is
immediate to show that $\varepsilon_{U}^{\ast}=0.5$. In fact, for any
$\varepsilon<0.5$, choosing $K=M$ we get that (\ref{epsilonun}) holds. This
proves that $\varepsilon_{U}^{\ast}\geq0.5$ and therefore $\varepsilon
_{U}^{\ast}=0.5.$ \

The following Theorem gives a lower bound for the FSBP of the estimate
$\widehat{\mu}_{n}$ defined in (\ref{OUREST}).

\begin{theorem}
\label{BDP}Let $\mathbf{W}$ be given by (\ref{Wcom}) and let $\mathbf{Z}%
=\{(\mathbf{x}_{i},y_{i}):i\in A\}.$ Suppose that $\widehat{\mathbf{\beta}%
}_{n}$=$\widetilde{\mathbf{\beta}}_{m}$ $(\mathbf{Z}),$ where $\widetilde
{\mathbf{\beta}}_{m}$ is a regression estimate for samples of size $m.$ Let
$\varepsilon_{1}>0$ be the FSBP at $\mathbf{Z}$\ of $\widetilde{\mathbf{\beta
}}_{m}$ and call$\ \varepsilon_{2}>0$ the UABP of $T_{L}.$\ Then the\ FSBP
\ $\varepsilon^{\ast}$ of the estimate $\widehat{\mu}_{n}$ \ at $\ \mathbf{W}$
satisfies
\[
\varepsilon^{\ast}\geq\varepsilon_{3}=\min(\varepsilon_{1},1-\sqrt
{1-\varepsilon_{2}}).
\]
\
\end{theorem}

In the next Section we introduce MM estimates of regression and
location. The maximum value of $\varepsilon_{1}$ for an MM estimate
of regression is $(n-c(G_{n}^{\ast}))/(2n)$, where $c(G)$ is defined
be (\ref{cG}) (see Martin et al. \cite{MaronnaMartin and Yohai}). In
Theorem \ref{locationbound} we show that maximum value of
$\varepsilon_{2}$ for an MM estimate of location is $0.5$. Then, if
$\ c(G_{n}^{\ast})/n$ is small, we can have have $\varepsilon_{3}$
close to $1-\sqrt{0.5}=0.293.$ A similar statement holds\ when we
change $T_{L}$ by the median.

\section{MM Regression and Location Functionals\label{SMM}}

Several robust estimates for the parameters of the regression model
(\ref{modelogral}) based on complete data $(\mathbf{x}_{1}%
,y_{1}),\dots,(\mathbf{x}_{n},y_{n})$ have been proposed. \ In this
paper we will consider MM estimates. These estimates were introduced
by Yohai \cite{Yohai87} for the linear model while  Fasano, Maronna,
Sued and Yohai \cite{Fasano et al}  extended these estimates   to
the case of nonlinear regression. For linear regression, MM
estimates may combine the highest possible breakdown point with an
arbitrarily high efficiency in \ the case of Gaussian errors. It
will be convenient \ to present MM-estimates of $\
$\textbf{$\beta$}$_{0}$ in their functional form, i.e., as \ a
functional $\mathbf{T}_{MM,\beta}(G)$ defined on \ a set \ of
distributions in $\ \mathbb{R}^{p+1}$, taking values in
$\mathbb{R}^{q}.$\ Given a sample $(\mathbf{x}_{1},y_{1}),\dots,(\mathbf{x}%
_{n},y_{n})$ \ the corresponding \ estimate of \textbf{$\beta$}$_{0}$ is given
by $\widehat{\mathbf{\beta}}_{MM}=\mathbf{T}_{MM,\beta}(G_{n}),$ where $G_{n}$
is the empirical distribution of the sample. As we explained in the
Introduction, we have excluded the intercept \ in \ model (\ref{modelogral}).
However in order to guarantee the consistency of the estimates without
requiring symmetric errors it is convenient to estimate an additional
parameter which can be naturally interpreted as an intercept or a center of
the error distribution. For this purpose\ put \textbf{$\xi$}$=($%
\textbf{$\beta$}$,\alpha)$ with $\alpha\in{\mathbb{R}}$, and define
$\underline{g}(\mathbf{x},${$\mathbf{\xi}$}$)=g(\mathbf{x},$\textbf{$\beta$%
}$)+\alpha$.

To define a \ regression MM functional $\mathbf{T}_{MM}(G)=(\mathbf{T}%
_{MM,\beta}(G),T_{MM,\alpha}(G))$\ two loss functions, $\ \rho_{0}$ $\ $\ and
$\rho_{1}$ are required$.$ The function $\rho_{0}$ is used \ to define \ a
dispersion functional $S(G)$ of the error distribution$.$Then $\mathbf{T}%
_{MM}$\ is defined as a regression M functional \ with loss function $\rho
_{1}$ and scale \ given by\ $\ S(G).$

Throughout this work, a \emph{bounded }$\rho$--\emph{function} is a function
$\rho\left(  t\right)  $ that is a continuous nondecreasing function of $|t|,$
such that $\rho(0)=0,$ $\rho\left(  \infty\right)  =1,$ and $\rho\left(
v\right)  <1$ implies that $\rho\left(  u\right)  <\rho(v)$ for $|u|<|v|.$ We
also assume that $\rho_{1}(t)\leq\rho_{0}(t)$ for all $t.$

We start by defining the dispersion functional. For any $\ $distribution $G$
of $(\mathbf{x},y)$ and $\mathbf{\xi}=($\textbf{$\beta$}$,\alpha),$ let
$S^{\ast}(G,\mathbf{\xi})$ \ be defined by%
\begin{equation}
\mathrm{E}_{G}\rho_{0}\left(  \frac{y-\underline{g}(\mathbf{x},\mathbf{\xi}%
)}{S^{\ast}(G,\mathbf{\xi})}\right)  =\delta, \label{Sestrella}%
\end{equation}
where $\delta\in(0,1)$. \ Then the dispersion functional $S(G)$ is defined by
\begin{equation}
S(G)=\min\limits_{\mathbf{\xi}\in B\times\mathbb{R}}S^{\ast}(G,\mathbf{\xi})
\label{SS}%
\end{equation}
and the MM estimating functional $\mathbf{T}_{\mathrm{MM}}(G)=(\mathbf{T}%
_{\mathrm{MM},\mathbf{\beta}}(G),T_{\mathrm{MM},\alpha}(G))$ by%
\begin{equation}
\mathbf{T}_{\mathrm{MM}}(G)=\arg\min\limits_{\mathbf{\xi}\in B\times
{\mathbb{R}}}\mathrm{E}_{G}\left[  \rho_{1}\left(  \frac{y-\underline
{g}(\mathbf{x},\mathbf{\xi})}{{S}(G)}\right)  \right]  . \label{TMM}%
\end{equation}

We can also consider another \ regression \ functional \ $\mathbf{T}%
_{\mathrm{S}}(G)=(\mathbf{T}_{\mathrm{S},\mathbf{\beta}}(G),T_{\mathrm{S}%
,\alpha}(G))$, \ called \emph{regression S functional}, as follows: \
\begin{equation}
\mathbf{T}_{\mathrm{S}}(G)=\arg\min\limits_{\mathbf{\xi}\in B\times
{\mathbb{R}}}\mathrm{E}_{G}\left[  \rho_{0}\left(  \frac{y-\underline
{g}(\mathbf{x},\mathbf{\xi})}{{S}(G)}\right)  \right]  . \label{TSS}%
\end{equation}

In the case of linear regression, the asymptotic breakdown point of both
$\mathbf{T}_{MM}$ and $\mathbf{T}_{S}$ is given by
\begin{equation}
\varepsilon^{\ast}=\min(\delta,1-\delta-c(G)), \label{BDPa}%
\end{equation}
where%
\begin{equation}
\label{cG}c(G)=\sup_{\mathbf{\gamma}\neq0,\mathbf{\gamma}\in\mathbb{R}^{p+1}%
}\mathrm{P}_{G}(\mathbf{\gamma}^{\prime}(\mathbf{x}^{\prime},1)^{\prime}=0).
\end{equation}

The maximum breakdown point \ occurs when $\delta=(1-c(G))/2$ and
its value is ($1-c(G))/2.$ It can be proved that this is the maximum
possible breakdown point for equivariant regression functionals. In
the case of nonlinear regression both $\mathbf{T}_{MM}$ and
$\mathbf{T}_{S}$ $\ $\ have also the same breakdown point but \ it
is not given by a simple closed expression (see Fasano
\cite{FasanoTesis}). \

Yohai \cite{Yohai87} showed that MM estimates for linear regression
may combine the highest possible breakdown point \ $(1-c(G))/2$ with
a Gaussian efficiency as high as desired. Instead, H\"{o}ssjer
\cite{Hossjer} showed that this is not \ possible for S estimates.
The maximum asymptotic Gaussian efficiency\ of an S estimate with
$\varepsilon^{\ast}=(1-c(G))/2$ is 0.33.

Let $(\mathbf{x},y)$ and $\ u$ satisfy model (\ref{modelogral}).\ Let
$\{\mathrm{G}_{n}^{\ast}\}$ be the sequence of empirical distribution
associated with observed pairs $(\mathbf{x}_{i},y_{i})$, i.e., those pairs
such that $\ a_{i}=1:$%
\begin{equation}
G_{n}^{\ast}=\frac{1}{\sum_{i=1}^{n}a_{i}}\,\sum_{i=i}^{n}a_{i}\delta
_{(\mathbf{x}_{i},y_{i})}. \label{pseudoempiricas}%
\end{equation}
Then we can estimate \textbf{$\beta$}$_{0}$ by
\begin{equation}
\widehat{\mathbf{\beta}}_{n}=\mathbf{T}_{MM,\mathbf{\beta}}(G_{n}^{\ast}).
\label{OURMM}%
\end{equation}

We can also choose as \ location functional $T_{L}$, whose value at $\mu
_{0}=T_{L}(F_{0})$ we want to estimate, a location MM functional. MM and S
location functionals are defined similarly to the regression case. Let
$\rho_{1}^{L}$ and $\rho_{0}^{L}$ be \ bounded$\ \rho$-functions. We start by
defining the dispersion functional. For any $\ $distribution $F$ of $y$ and
$\mu\in\mathbb{R}\ $ let $S_{L}^{\ast}(F,\mu)$ \ be defined by%
\[
\mathrm{E}_{F}\rho_{0}^{L}\left(  \frac{y-\mu}{S_{L}^{\ast}(F,\mathbf{\xi}%
)}\right)  =\delta,
\]
where $\delta\in(0,1)$. \ Then the dispersion functional $S_{L}(F)$ is defined
by
\[
S_{L}(F)=\min\limits_{\mu\in\mathbb{R}}S_{L}^{\ast}(F,\mu)
\]
and \ the MM location functional $T_{\mathrm{MM}}^{L}(F)$ by%
\begin{equation}
T_{\mathrm{MM}}^{L}(F)=\arg\min\limits_{\mu\in\mathbb{R}}\mathrm{E}_{F}\left[
\rho_{1}^{L}\left(  \frac{y-\mu}{{S}_{L}(F)}\right)  \right]  .
\label{del loc}%
\end{equation}

The \ S location functional $T_{\mathrm{S}}^{L}(F)$\ \ is defined similarly to
the regression S functional. We denote by $\mu_{00}=T_{\mathrm{S}}^{L}(F_{0})$
and $\mu_{01}=T_{\mathrm{MM}}^{L}(F_{0})$, whenever they are well defined.
Location MM estimates may also combine \ high breakdown point with high
Gaussian efficiency and their breakdown point is given by $\varepsilon^{\ast
}=\min(\delta,1-\delta).$

For the validity \ of assumptions A1-A4, the $\rho$-functions used to define
the location and regression MM functionals should satisfy assumptions R1 and
R2 below.

\textbf{R1} For some $m,$ $\rho(u)=1$ iff $|u|\geq m,$ and $\log(1-\rho)$ is
concave on $(-m,m)$ .

\textbf{R2} $\rho$ is twice continuously differentiable

A family of very popular bounded $\rho-$function satisfying R0,R1 and R2 is
Tukey's bisquare family:
\begin{equation}
\rho_{k}^{T}\left(  u\right)  =1-\left(  1-\left(  \frac{u}{k}\right)
^{2}\right)  ^{3}I(|u|\leq k)
\end{equation}

for $\ k>0.$

We denote by $\psi_{0}$, $\psi_{1},\psi_{0}^{L}$ and $\psi_{1}^{L}$ the
derivatives of $\rho_{0}$, $\rho_{1},\rho_{0}^{L}$ and $\rho_{1}^{L}$.\ \ Put
$\alpha_{01}=T_{MM,\alpha}(G_{0}^{\ast}),$ $\alpha_{00}=T_{S,\alpha}%
(G_{0}^{\ast})$ and $\sigma_{0}=S(G_{0}^{\ast})$

Both regression and location MM and S functionals are studied in
detail in Fasano et al. \cite{Fasano et al}. There we can find
sufficient conditions for weak continuity and Fisher-consistency.\
Moreover, a weak differentiability notion involving the influence
function of the functionals is also developed. \ This notion allows
to obtain asymptotic expansions, like those required in
(\ref{TRexpansion}) and (\ref{TLexpansion}). The following numbers
will be used to derive the influence functions of the regression
functionals:
\[
a_{0i}=\mathrm{E}_{G_{0}^{\ast}}\psi_{i}^{\prime}\left(  (y-g(\mathbf{x}%
,\mathbf{\beta}_{0})-\alpha_{0i})/\sigma_{0}\right)  =\mathrm{E}_{K_{0}}%
\psi_{i}^{\prime}\left(  (u-\alpha_{0i})/\sigma_{0}\right)  ,i=0,1,
\]%
\[
e_{0i}=\mathrm{E}_{K_{0}}\left[  \psi_{i}^{\prime}\left(  (u-\alpha
_{0i})/\sigma_{0}\right)  (u-\alpha_{0i})/\sigma_{0}\right]  ,i=0,1,
\]%
\[
d_{0}=\mathrm{E}_{K_{0}}\left[  \psi_{0}\left(  (u-\alpha_{00})/\sigma
_{0}\right)  (u-\alpha_{00})/\sigma_{0}\right]  \quad\hbox{and}\quad
\mathbf{b}_{0}=\mathrm{E}_{G_{0}^{\ast}}\dot{g}(\mathbf{x},\mathbf{\beta}%
_{0}).
\]
Similarly we define $a_{0i}^{L}$, $\ e_{0i}^{L}$, $d_{0}^{L}$ and $\sigma
_{0}^{L}$ replacing $\psi_{_{i}}$ by $\psi_{i}^{L}$, $K_{0}$ by $F_{0}$,
$g(\mathbf{x},$\textbf{$\beta$}$_{0})$ by $0$, $\alpha_{0i}$ by $\mu_{0i}$ and
$\sigma_{0}$ by $\sigma_{0}^{L}=S_{L}(F_{0})$. We denote by $A_{0}$ the
covariance matrix of $\dot{g}(\mathbf{x},$\textbf{$\beta$}$_{0})$ under
$Q_{0}^{\ast}.$

Theorems \ref{asslreg} and \ref{assloc} summarize the results for MM
functionals of regression and location, respectively.

\begin{theorem}
\label{asslreg}Let $\rho_{0}$ and $\rho_{1}$ be bounded $\rho$-functions
satisfying $\ $R1, with $\rho_{1}\leq\rho_{0}$. Assume that $K_{0}$ has a
strongly unimodal density and that (\ref{IDCOND}) holds replacing $Q_{0}$ by
$Q_{0}^{\ast}$. We will consider that either (a) $B$ is compact \ or (b)
\ $g(\mathbf{x},$\textbf{$\beta$}$)=\mathbf{\beta}^{\prime}\mathbf{x}$ and
$\delta<1-c(G_{0}^{\ast})$. Then

\begin{itemize}
\item[(i)] $\lim_{n\rightarrow\infty}$\ $\mathbf{T}_{\mathrm{MM,}%
\mathbf{\beta}}(G_{n}^{\ast})=$\textbf{$\beta$}$_{0}$ $\ $\ a.s. and therefore
A1 is satisfied.

\item[(ii)] Assume also that $a_{00},$\ $a_{01}$ and $d_{0}$ are different
from $0$, that A0 holds and that $\rho_{0}$ and $\rho_{1}$ satisfies R2. Then
(\ref{TRexpansion}) holds \ with $I_{R}(\mathbf{x},y)=I_{\mathbf{T}%
_{MM,\mathbf{\beta}},G_{0}^{\ast}}(\mathbf{x,}y\mathbf{)/}\mathrm{E}(a_{1})$ ,
where $I_{\mathbf{T}_{MM,\mathbf{\beta}},G_{0}^{\ast}}(\mathbf{x,}y\mathbf{)}
$ is the influence function of $\mathbf{T}_{MM,\mathbf{\beta}}$ at
$G_{0}^{\ast}.$ Moreover, we have that
\begin{equation}
I_{\mathbf{T}_{MM,\mathbf{\beta}},G_{0}^{\ast}}(\mathbf{x,}y\mathbf{)}%
=\frac{\sigma_{0}}{a_{01}}\psi_{1}\left(  \frac{y-\underline{g}(\mathbf{x}%
,(\mathbf{\beta}_{0},\alpha_{01}))}{\sigma_{0}}\right)  A_{0}^{-1}(\dot
{g}(\mathbf{x},\mathbf{\beta}_{0})-\mathbf{b}_{0}), \label{ICEX2}%
\end{equation}
and therefore A2 holds.\bigskip
\end{itemize}
\end{theorem}

\begin{theorem}
\label{assloc} Let $\rho_{0}^{L}$ and $\rho_{1}^{L}$ be bounded $\rho
$-functions satisfying $\ $R1, with $\rho_{1}^{L}\leq\rho_{0}^{L}$. Assume
that $\ F_{0}$ has a strongly unimodal density. Then

\begin{itemize}
\item[(i)] There is only one value $\mu_{01}=T_{MM}^{L}$ $(F_{0})$ that
attains the minimum \ at (\ref{del loc}), $T^{L}_{MM}$ is continuous at
$F_{0}$, and so A3 holds. In the case that \ $F_{0}$ is symmetric around
$\nu_{0},$ we have $\mu_{01}=\nu_{0}$.

\item[(ii)] Assume also A0, that $\rho_{0}^{L}$ and $\rho_{1}^{L}$ satisfy
R2\ and \ that $a_{00}^{L},$\ $a_{01}^{L}$ and $d_{0}^{L}$ are different from
$0.$ Then (\ref{TLexpansion}) holds \ when $I_{L}(y)$ is the influence
function of $T_{MM}^{L}$ at $F_{0}.$ Moreover we have
\begin{equation}
I_{L}(y)=\frac{\sigma_{0}^{L}}{a_{01}^{:L}}\psi_{1}^{L}\left(  \frac
{y-\mu_{01})}{\sigma_{0}^{L}}\right)  -\frac{e_{01}^{L}\ \sigma_{0}^{L}%
}{a_{01}^{L}d_{0}^{L}}\left(  \rho_{0}^{L}\left(  \frac{y-\mu_{00}}{\sigma
_{0}^{L}}\right)  -\delta\right)  , \label{ifloc1}%
\end{equation}
and therefore A4 holds.

\item[(iii)] \ In case \ that $F_{0}$ is symmetric with respect to $\nu_{0}$
we have $e_{0}=0$ and
\[
I_{L}(y)=\frac{\sigma_{0}^{L}}{a_{01}^{:L}}\psi_{1}^{L}\left(  \frac{y-\nu
_{0})}{\sigma_{0}^{L}}\right)  .
\]

\end{itemize}
\end{theorem}

To end this Section, we state the announced result regarding the uniform bound
required for the location functional in order to deduce a lower bound for the
FSBD of $\widehat\mu_{n}$, introduced in Section \ref{secBDP}.

\begin{theorem}
\label{locationbound} Let $T_{MM}^{L}$ be an MM location functional. Then its
uniform asymptotic breakdown point is $\varepsilon_{U}^{\ast}=\min
(1-\delta,\delta)$.
\end{theorem}

\section{Monte Carlo study \label{SECMC}}

In order to assess how the proposed robust method compares to the classical
procedure that uses as $\widehat{\mathbf{\beta}}_{n}$ the least squares and as
$T_{L}$ the mean functional, we performed a Monte Carlo study. We consider the
following model
\[
y_{i}=3x_{i1}+...+3x_{i5}+u_{i},1\leq i\leq100,
\]
where $x_{i1},...,x_{i5}$ are i.i.d. random variables with uniform
distribution in the interval $[0,1]$, $u_{i}$ are standardized normal
variables ($u_{i}\sim{\mathcal{N}}(0,1)$) and $\beta_{1}=\beta_{2}%
=...=\beta_{5}=3$. The missingness indicators $a_{i}$ were generated using a
logistic model. Let $\mathbf{x}_{i}=(x_{i1},...,x_{i5})$, then
\[
\log\frac{\mathrm{P}(a_{i}=1|\mathbf{x}_{i})}{1-\mathrm{P}(a_{i}%
=1|\mathbf{x}_{i})}=0.57(x_{i1}+...+x_{i5}).
\]
Using this model and the distribution of the covariables, we have
$\mathrm{P}(a_{i}=1)=0.80$.

We study (a) the case with no outlier contamination and (b) the case where
10\% of the observations $(\mathbf{x,}$ $y_{i})$'s with $a_{i}=1$ are replaced
by $(\mathbf{x}^{\ast},y^{\ast})$, with $\mathbf{x}^{\ast}=(x^{\ast
},...,x^{\ast})$. We take two values for $x^{\ast}$: $1$ and $3$, and for
$y^{\ast}$ we take a grid of values over the interval $[8,50]$, with steps of
$0.20$. For each case we performed 1000 replications. We consider four
functionals $T_{L}$ : (i) the mean (MEAN in Figure 1), (ii) the median (MEDIAN
in Figure 1) (iii) an MM location functional \ with $\rho_{i}^{L}%
=\rho_{T,k_{i}},$ $k_{0}$=$1.57$, $k_{1}=3.88$ and $\delta=0.5$ $.$ The
corresponding location estimate has a Gaussian asymptotic efficiency of 90\%
\ (MM90 in Figure 1). (iv) Finally\ we study an MM location functional defined
as in (iii) \ with constants \ $k_{0}$=$1.57$, $k_{1}=4.68$ and $\delta=0.5$.
This location estimate \ has a Gaussian asymptotic efficiency of 95\% (MM95 in
Figure 1). Note that \ in the case in which there is no outlier contamination,
the distribution $F_{0}$ is symmetric with center of symmetry $7.5$, and then
$T_{L}(F_{0})=\mathrm{E}(y)=7.5$ in the four cases. When $T_{L}$ is the mean,
$\widehat{\mathbf{\beta}}_{n}$ is the least squares (LS) estimate. In the
other 3 cases $\widehat{\mathbf{\beta}}_{n}$ is an MM estimate with $\rho
_{i}=\rho_{T,k_{i}},$ $k_{0}$=$1.57$, $k_{1}=3.44$ and $\delta=0.5.$ This
estimate has an asymptotic efficiency of 85\% in the case of Gaussian errors
and breakdown point close to $0.5.$ In Table 1 we show the mean square errors
(MSE), and the relative efficiencies of the four estimates when there is no
outlier contamination. In Figure 1 we plot the MSE of the four estimates
\ under outlier contamination. \bigskip

\begin{center}
Table 1. MSE and efficiencies without outliers\bigskip
\end{center}

\[%
\begin{tabular}
[c]{lllll}%
Estimates & MEAN & MEDIAN & MM90 & MM95\\
MSE & 0.047 & 0.056 & 0.051 & 0.049\\
Efficiency & 100\% & 83\% & 91\% & 95\%
\end{tabular}
\
\]
%


\begin{figure}
[ptb]
\begin{center}
\includegraphics[
height=4.5766in,
width=5.1171in
]%
{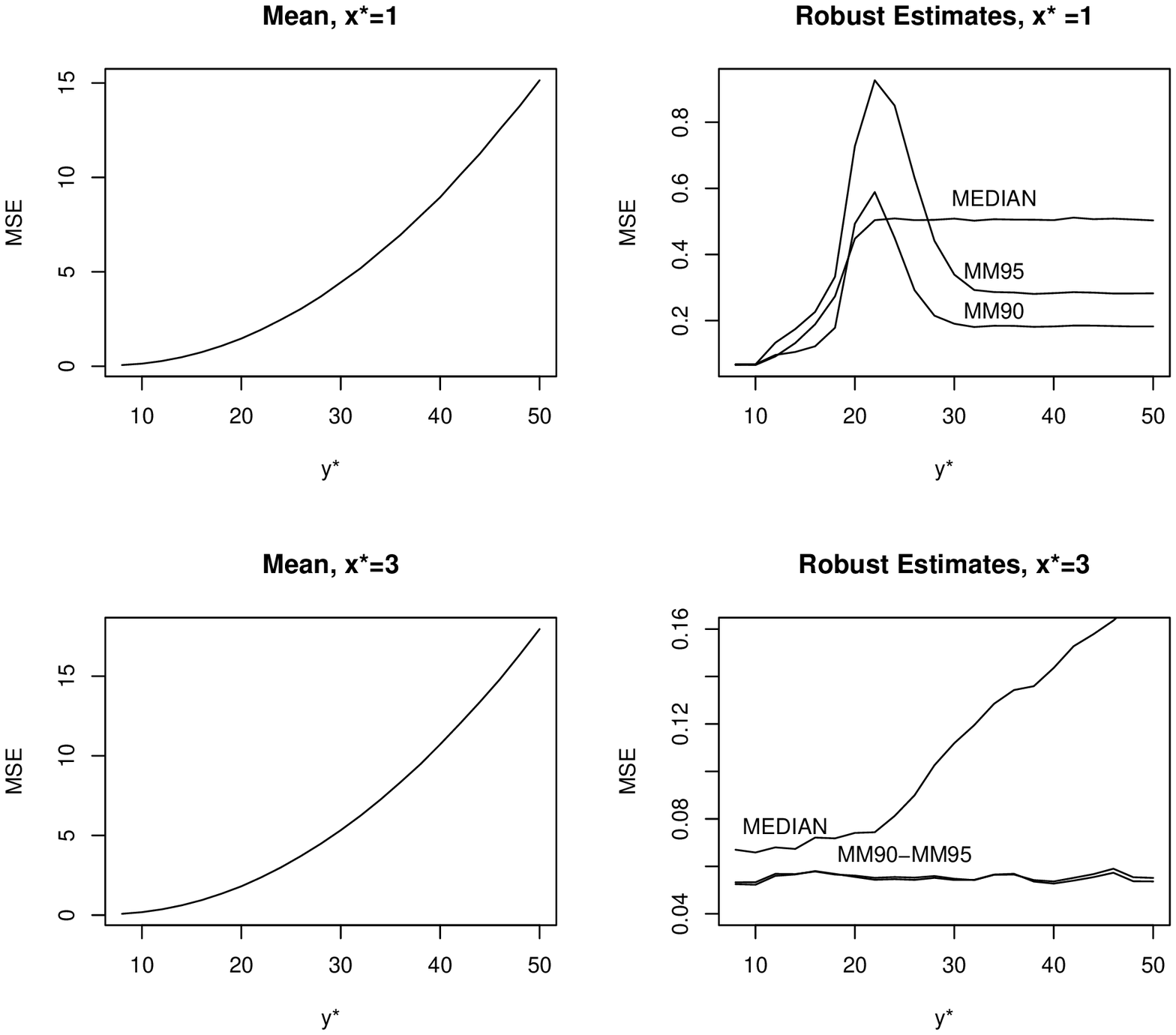}%
\end{center}
\end{figure}
\bigskip As expected, when there are no outliers the classical estimate based
on the mean is the most efficient, but the estimates based on the MM estimates
are highly efficient too. The estimate based on the median is less efficient,
but its efficiency is larger than that of the \ sample median which is 64\%.
Note that the estimate based on the median is an U-statistics similar to the
Hodges--Lehmann estimate, which is also more efficient than the median.

When there are outliers, we observe that the MSE of the estimate based on the
mean increases beyond any limit, while for the robust estimates the MSE
remains bounded. In the case of $x^{\ast}=1$ the MSE of MM95 is larger than
those of MEDIAN and MM90. For $x^{\ast}=3$ the MSE of MEDIAN is larger than
those of the other two robust estimates. The MSEs of MM90 and MM95 are
practically the same. Based on these results we recommend to use MM90 which
has a very good behavior with and without outliers.

\section*{ Acknowledgements} The authors would like to thank Graciela
Boente and Sara van der Geer for valuable discussions and
suggestions and also to Damian Scherlis for his careful reading of
the first manuscript.

\section{Appendix{\label{Apendice}}}

The following result plays a crucial role in the proof of Theorem
\ref{lemaconvdebil1}.

\begin{lemma}
\label{convdebil} Let $\{\mathbf{z}_{i}\}$ be a sequence of i.i.d. random
vectors taking values in ${\mathbb{R}}^{k}$ and let $h:{\mathbb{R}}^{k}%
\times{\mathbb{R}}^{q}\rightarrow{\mathbb{R}}$ be a continuous function.
Assume that $\widehat{\mathbf{\beta}}_{n}$ is a strongly consistent sequence
of estimators of $\beta_{0}\in{\mathbb{R}}^{q}$. Denote by $\widehat{H}_{n}$
the empirical distribution at $h(\mathbf{z}_{i},\widehat{\mathbf{\beta}}_{n}%
)$, $1\leq i\leq n$ and by $H_{0}$ the distribution of $h(\mathbf{z}_{1}%
,\beta_{0})$. Then $\widehat{H}_{n}$ converges weakly to $H_{0}$ a.s., i.e.%
\begin{equation}
\mathrm{P}(\widehat{H}_{n}\rightarrow_{w}H_{0})=1. \label{convdebilae}%
\end{equation}

\end{lemma}

\begin{proof} Recall that weak convergence is characterized by the following
property:
\[
H_{n}\rightarrow H\hbox{weakly}\Leftrightarrow\int f\,dH_{n}\rightarrow\int
f\,dH,\quad\forall f\in{\mathcal{C}}_{B}({\mathbb{R}}),
\]
where ${\mathcal{C}}_{B}({\mathbb{R}})$ denotes the set of continuous bounded
functions. Denote by $\tilde{H_{n}}$ the empirical distribution at
$h(\mathbf{z}_{i},\beta_{0})$, for $1\leq i\leq n$. By the Glivenko-Cantelli
Theorem, $\tilde{H_{n}}$ converges uniformly to $H_{0}$, a.s. and so \ it also
converges weakly a.s. Then, it remains to find a set of probability one where
\[
\lim\limits_{n\rightarrow\infty}\left\vert \int f\,d\widehat{H}_{n}-\int
f\,d\tilde{H_{n}}\right\vert =0,\text{ }\forall f\in{\mathcal{C}}%
_{B}({\mathbb{R}}).
\]
Observe that
\[
\int f\,d\widehat{H}_{n}=\frac{1}{n}\sum_{i=1}^{n}f\left(  h(\mathbf{z}%
_{i},\widehat{\beta}_{n})\right)  ,\text{ }\int f\,d\tilde{H_{n}}=\frac{1}%
{n}\sum_{i=1}^{n}f\left(  h(\mathbf{z}_{i},\widehat{\beta}_{0})\right)  ,
\]
and so
\[
\left\vert \int f\,d\widehat{H}_{n}-\int f\,d\tilde{H_{n}}\right\vert
\leq\frac{1}{n}\sum_{i=1}^{n}\left\vert f\left(  h(\mathbf{z}_{i}%
,\widehat{\beta}_{n})\right)  -f\left(  h(\mathbf{z}_{i},\widehat{\beta}%
_{0})\right)  \right\vert I_{\{|\mathbf{z}_{i}|\leq K\}}+2||f||_{\infty}%
\frac{1}{n}\sum_{i=1}^{n}I_{\{|\mathbf{z}_{i}|>K\}}.
\]

Put $C_{K}=\{(\mathbf{z},\beta):||\mathbf{z}||\leq K,||\beta-\beta_{0}%
||\leq1\}$. We have that $f\circ h:C_{K}\rightarrow\mathbb{R}$ is uniformly
continuous and so, given $\varepsilon>0$, there exists $\delta>0$ such that if
$(\mathbf{z}_{i},\beta_{i})\in C_{K}$ and $||(\mathbf{z}_{1},\beta
_{1})-(\mathbf{z}_{2},\beta_{2})||\leq\delta$, then $|f(h(\mathbf{z}_{1}%
,\beta_{1}))-f(h(\mathbf{z}_{1},\beta_{2}))|\leq\varepsilon$. With probability
one there exists a \ random integer $n_{0}$ $\ $\ such that $|\widehat
{\mathbf{\beta}}_{n}-$\textbf{$\beta$}$_{0}|\leq\delta$ for all $n\geq n_{0}$.
Then we get
\[
\left\vert \int f\,d\widehat{H}_{n}-\int f\,d\tilde{H_{n}}\right\vert
\leq\varepsilon+2||f||_{\infty}\frac{1}{n}\sum_{i=1}^{n}I_{\{|\mathbf{z}%
_{i}|>K\}},
\]
for all $n\geq n_{0}$. Assume also that
\[
\frac{1}{n}\sum_{i=1}^{n}I_{\{|\mathbf{z}_{i}|>K\}}\rightarrow\mathrm{P}%
\left(  |\mathbf{z}_{1}|>K\right)  ,\forall K.
\]
Then, with probability one
\[
\lim\limits_{n\rightarrow\infty}\left\vert \int f\,d\widehat{H}_{n}-\int
f\,d\tilde{H_{n}}\right\vert \leq\varepsilon+2||f||_{\infty}\,\mathrm{P}%
\left(  |\mathbf{z}_{1}|>K\right)  ,\forall\varepsilon>0\,,\forall K.
\]
To get the desired result, let $\varepsilon\rightarrow0$ and $K\rightarrow
\infty$.

\end{proof}

The following results will be used throughout the proofs of the Theorems
stated in the previous Sections. We start proving that the convolution
preserves weak continuity.

\begin{lemma}
\label{convol} Assume that $K_{n}\rightarrow_{w}K_{0}$ and $\ R_{n}%
\rightarrow_{w}R_{0}.$ Then $K_{n}\ast R_{n}\rightarrow_{w}K_{0}\ast R_{0}.$
\end{lemma}

\begin{proof} Let $(U,V)$ be independent random variables, both with uniform
distribution on $[0,1]$. Given a distribution function $F$, denote
by $F^{-1}$ the generalized inverse function of $F$, whose value at
$t$ is given by the infimum of the set $\{s:t\leq F(s)\}$. Consider
$U_{n}=K_{n}^{-1}(U)$ and $V_{n}=R_{n}^{-1}(V)$. It is known that
(i) $U_{n}$ and $V_{n}$ are distributed according $K_{n}$ and
$R_{n}$, respectively and (ii) $U_{n}$ and $V_{n}$ converge a.s. to
$U_{0}=K_{0}^{-1}(U)$ and $V_{0}=R_{0}^{-1}(V)$, respectively (see
Theorem 25.6 (Billingsley (1995)) for details). Then $U_{n}+V_{n}$
converges a.s. to $U_{0}+V_{0}$, and then the convergence holds also
in distribution. The independence between $U$ and $V$ implies that
$U_{n}+V_{n}\sim K_{n}\ast R_{n}$ while $U_{0}+V_{0}\sim K_{0}\ast
R_{0}$, proving the Lemma.\end{proof}

\begin{lemma}
\label{g-c}Consider $\{(a_{i},\mathbf{z}_{i})\}$ i.i.d. random vectors, with
Bernoulli $a_{i}$ \ and $\mathbf{z}_{i}\in\mathbb{R}^{h}$ Then
\begin{equation}
\sup_{z\in{\mathbb{R}}^{h}}\left\vert \frac{1}{n}\sum_{i=1}^{n}a_{i}%
I_{\{\mathbf{z}_{i}\leq\mathbf{z\}}}-\mathrm{E}\left[  a_{1}I_{\{\mathbf{z}%
_{1}\leq z\}}\right]  \right\vert =0,\;\text{as..}%
\end{equation}

\end{lemma}

\begin{proof} \ Note that%
\begin{equation}
a_{i}I_{\{\mathbf{z}_{i}\leq\mathbf{z\}}}=I_{\{\mathbf{z}_{i}\leq\mathbf{z\}}%
}-I_{\{\mathbf{z}_{i}\leq\mathbf{z,}a_{i}\leq0\mathbf{\}}}. \label{GC0}%
\end{equation}
\ By \ the Glivenko-Cantelli Theorem \ we have%
\begin{equation}
\sup_{z\in{\mathbb{R}}^{h}}\left\vert \frac{1}{n}\sum_{i=1}^{n}I_{\{\mathbf{z}%
_{i}\leq\mathbf{z\}}}-\mathrm{E}\left[  I_{\{\mathbf{z}_{1}\leq z\}}\right]
\right\vert =0,\;\text{a.s.} \label{GC1}%
\end{equation}
and%
\begin{equation}
\sup_{z\in{\mathbb{R}}^{h}}\left\vert \frac{1}{n}\sum_{i=1}^{n}I_{\{\mathbf{z}%
_{i}\leq\mathbf{z,}a_{i}\leq0\mathbf{\}}}-\mathrm{E}\left[  I_{\{\mathbf{z}%
_{1}\leq\mathbf{z,}a_{1}\leq0\mathbf{\}}}\right]  \right\vert
=0,\;\text{a.s.} \label{GC2}%
\end{equation}
From (\ref{GC0}),(\ref{GC1}) and (\ref{GC2}) we get%
\[
\sup_{z\in{\mathbb{R}}^{h}}\left\vert \frac{1}{n}\sum_{i=1}^{n}a_{i}%
I_{\{\mathbf{z}_{i}\leq\mathbf{z\}}}-\mathrm{E}\left[  I_{\{\mathbf{z}_{1}\leq
z\}}-I_{\{\mathbf{z}_{1}\leq\mathbf{z,}a_{1}\leq0\mathbf{\}}}\right]
\right\vert .
\]
and by \ applying (\ref{GC0}) \ to $i=1$ the Lemma follows.
\end{proof}
The proof of the following Lemma is similar to that of Lemma 4.2
presented by  Yohai in \cite{Yohai85}. It suffices to replace the
law of large numbers \ for i.i.d., variables by the same law for U
statistics.

\begin{lemma}
\label{victor4.2} Assume that $\{\mathbf{z}_{i}\}$ are i.i.d. random vectors
taking values in ${\mathbb{R}}^{k}$, with common distribution $Q$. Let
$f:{\mathbb{R}}^{k}\times{\mathbb{R}}^{k}\times{\mathbb{R}}^{h}\rightarrow
{\mathbb{R}}$ be a continuous function. Assume that for some $\delta>0$ we
have that%
\[
\mathrm{E}\sup_{||\lambda-\lambda_{0}||\leq\delta}|f(\mathbf{z}_{1}%
,\mathbf{z}_{2},\lambda)|<\infty
\]
$\ $and \ that $\widehat{\lambda}_{n}\rightarrow\lambda_{0}$ a.s. Then\
\begin{equation}
\frac{1}{n^{2}}\sum_{j=1}^{n}\sum_{1=1}^{n}f(\mathbf{z}_{i},\mathbf{z}%
_{j},\widehat{\lambda}_{n})\rightarrow\mathrm{E}f(\mathbf{z}_{1}%
,\mathbf{z}_{2},\lambda_{0})\,\text{\ a.s.} \label{lgnu}%
\end{equation}
\bigskip
\end{lemma}

\textsc{Proof of Theorem \ref{lemaconvdebil1}. }According to Lemma
\ref{convol}, it only remains to prove the a.s. weak convergence of
$\widehat{R}_{n}$ $\ $\ and $\widehat{K}_{n}$ \ to $\ R_{0}$ and $K_{0}$
respectively.\ The a.s. weak convergence of $\widehat{R}_{n}$ to $R_{0}$
follows from Lemma \ref{convdebil}, putting $\mathbf{z}=\mathbf{x}$ and
$h(\mathbf{z},\beta)=g(\mathbf{x},\beta)$. Weak convergence of $(\widehat
{K}_{n})_{n\geq1}$ to $K_{0}$ requires an extra argument. If $\mathbf{z}%
=(\mathbf{x},y)$ and $h(\mathbf{z},\beta)=y-g(\mathbf{x},\beta)$, we get that
\[
\widehat{K}_{n}(u)=\frac{1}{\sum_{i=1}^{n}a_{i}}\,{\sum_{i=1}^{n}%
a_{i}\,I_{\{h(\mathbf{z}_{i},\widehat{\mathbf{\beta}}_{n}))\leq u\}}}.
\]
By Lemma \ref{g-c}, we obtain%
\[
\sup_{u\in{\mathbb{R}}}\left\vert \frac{1}{n}\sum_{i=1}^{n}a_{i}I_{\{u_{i}\leq
u\}}-\mathrm{E}\left[  a_{1}I_{\{u_{1}\leq u\}}\right]  \right\vert
=0\;\text{a.s.}%
\]
Since $a_{1}$ and $u_{1}$ are independent, we conclude that
\[
\sup_{u\in{\mathbb{R}}}\left\vert \frac{1}{\sum_{i=1}^{n}a_{i}}\sum_{i=1}%
^{n}a_{i}I_{\{u_{i}\leq u\}}-K_{0}(u)\right\vert =0\;\text{a.s.}%
\]
and $\ $\ then $\sum_{i=1}^{n}a_{i}I_{\{u_{i}\leq u\}}/\sum_{i=1}^{n}a_{i}$
converges weakly to $K_{0}$ a.s. An argument \ similar to the one used in
Lemma \ref{convdebil} \ shows that with probability one \ we have%
\[
\lim\limits_{n\rightarrow\infty}\left\vert \int f\,d\widehat{K}_{n}-\int
f\,d\tilde{K_{n}}\right\vert =0,\quad\forall f\in{\mathcal{C}}_{B}%
({\mathbb{R}}),
\]
proving \ the \ a.s. weak convergence of $\widehat{K}_{n}$ to $K_{0}$ . This
concludes the proof of \ part (a) of Theorem \ref{lemaconvdebil1}. (b) is an
immediate consequence of weak continuity of $T_{L}$. $\square$\bigskip

\textsc{Proof of Theorem \ref{asym}. }According to A4, we have that
\[
\sqrt{n}(\widehat{\mu}_{n}-\mu_{0})\;=\;\sqrt{n}\Big\{T_{L}(\widehat{F}%
_{n})-T_{L}(F_{0})\Big\}\;=\;\sqrt{n}\mathrm{E}_{\widehat{F}_{n}}%
I_{L}(y)\,+\,o_{P}(1).
\]
Note that
\[
\mathrm{E}_{\widehat{F}_{n}}I_{L}(y)=\frac{1}{\eta_{n}n^{2}}{\sum_{j=1}^{n}%
}{\sum_{i=1}^{n}}a_{i}I_{L}(y_{i}-g(\mathbf{x}_{i},\widehat{\mathbf{\beta}%
}_{n})+g(\mathbf{x}_{j},\widehat{\mathbf{\beta}}_{n})),
\]
where $\eta_{n}=\sum_{i=1}^{n}a_{i}/n$. Since $\eta_{n}\rightarrow
\mathrm{E}[a_{i}]=\eta$, to prove Theorem \ref{asym}, it is enough to show
that%
\[
V_{n}\rightarrow_{d}N(0,(\eta\tau)^{2}),
\]
where%

\[
V_{n}=\frac{1}{n^{3/2}}%
{\displaystyle\sum_{j=1}^{n}}
{\displaystyle\sum_{j=1}^{n}}
a_{i}I_{L}(y_{i}-g(\mathbf{x}_{i},\widehat{\mathbf{\beta}}_{n})+g(\mathbf{x}%
_{j},\widehat{\mathbf{\beta}}_{n})).
\]
Performing a Taylor expansion, we can write
\[
V_{n}=d_{n}+\mathbf{c}_{n}^{\prime}n^{1/2}(\widehat{\mathbf{\beta}}%
_{n}-\mathbf{\beta}_{0}),
\]
where
\[
d_{n}=\frac{1}{n^{3/2}}%
{\displaystyle\sum_{i=1}^{n}}
{\displaystyle\sum_{j=1}^{n}}
a_{i}I_{L}(u_{i}+g(\mathbf{\beta}_{0},\mathbf{x}_{j}))
\]
and%
\[
\mathbf{c}_{n}=\frac{1}{n^{2}}\sum_{i=1}^{n}\sum_{j=1}^{n}\ell(a_{i}%
,\mathbf{x}_{i},y_{i},a_{j},\mathbf{x}_{j},y_{j},\beta_{n}^{\ast})
\]
with \textbf{$\beta$}$_{n}^{\ast}$ between $\widehat{\mathbf{\beta}}_{n}$ and
\textbf{$\beta$}$_{0}$, and
\[
\ell(a_{i},\mathbf{x}_{i},y_{i},a_{j},\mathbf{x}_{j},y_{j},\beta)=a_{i}%
\,I_{L}^{\prime}(y_{i}-g({\mathbf{\beta}},\mathbf{x}_{i})+{g}({\mathbf{\beta}%
},\mathbf{x}_{j}))\,\{\dot{g}({\mathbf{\beta}},\mathbf{x}_{j})-\dot
{g}({\mathbf{\beta}},\mathbf{x}_{i})\}.
\]
By Lemma \ref{victor4.2}
\begin{equation}
\mathbf{c}_{n}\rightarrow\mathbf{c}=\mathrm{E}\ell(a_{1},\mathbf{x}_{1}%
,y_{1},a_{2},\mathbf{x}_{2},y_{2},\beta_{0})\text{ a.s.} \label{asf1}%
\end{equation}
\ From the U-statistics projection Theorem \ we get
\begin{equation}
d_{n}=\frac{1}{n^{1/2}}%
{\displaystyle\sum\limits_{i=1}^{n}}
e(\mathbf{x}_{i},u_{i},a_{i})+f(\mathbf{x}_{i})+o_{P}(1). \label{asf2}%
\end{equation}
Finally, using (\ref{TRexpansion}), we get that
\[
V_{n}=\frac{1}{n^{1/2}}%
{\displaystyle\sum\limits_{i=1}^{n}}
e(\mathbf{x}_{i},u_{i},a_{i})+f(\mathbf{x}_{i})+a_{i}{\mathbf{c}}^{{\prime}%
}{I_{R}(\mathbf{x}_{i},y_{i})}+o_{P}(1),
\]
and using the Central Limit Theorem we get (\ref{Vnasnor}). $\square$\bigskip

To prove Theorem \ref{median}\ \ we need \ an asymptotic expansion for
$n^{1/2}(\widehat{\mu}_{n}-\mu_{0})$. Let $\mathbf{z}_{i}=(a_{i}%
,\mathbf{x}_{i},y_{i})$ and consider
\begin{equation}
\Psi_{1}(\mathbf{z}_{i},\mathbf{z}_{j},\beta,\mu)=a_{i}\text{sign}\left(
g(\mathbf{x}_{j},\beta)+y_{i}-g(\mathbf{x}_{i},\beta)-\mu\right)  ,
\label{Psi}%
\end{equation}%
\[
\mathbf{\Lambda}_{1}(\beta,\mu)=\mathrm{E}\Psi(\mathbf{z}_{1},\mathbf{z}%
_{2},\beta,\mu)\;,\Lambda_{1\mathbf{\beta}}(\mathbf{\beta},\mu\mathbf{)}%
=\frac{\partial\Lambda_{1}(\mathbf{\beta},\mu)}{\partial\mathbf{\beta}},\text{
}\Lambda_{1\mu}(\mathbf{\beta},\mu\mathbf{)}=\frac{\partial\Lambda
_{1}(\mathbf{\beta},\mu)}{\partial\mu},
\]
and
\[
J_{n}(\beta,\mu)=\frac{1}{n^{3/2}}\sum_{j=1}^{n}\sum_{i=1}^{n}\Psi
_{1}(\mathbf{z}_{i},\mathbf{z}_{j},\beta,\mu).
\]
The independence between $a_{1}$ and $(u_{1},\mathbf{x}_{2})$ and the fact
that $u_{1}+g(\mathbf{x}_{2},\mathbf{\beta}_{0})$ has distribution $F_{0}$,
allow to conclude that \ $\Lambda_{1}(\beta_{0},\mu_{0})=\mathrm{E}%
(a_{1})\mathrm{E}_{F_{0}}\text{sign}(y-\mu_{0})=0$.\ Since $I_{T_{\text{med}%
},F_{0}}(y)$ is not differentiable we have to use \ an extra argument to
obtain an asymptotic linear expansion for $n^{1/2}(\widehat{\mu}_{n}-\mu_{0}%
)$. To this purpose, the following Lemma is crucial. It is related
to a very general linear expansion satisfied by empirical processes
based on U-statistics.
\begin{lemma}
\label{tayfa} Suppose the same assumptions \ as in Theorem \ref{median}. Then
if $n^{1/2}$ $(\overline{\beta}_{n}-\beta_{0})$ and $n^{1/2}(\overline{\mu
}_{n}-\mu_{0})$ are bounded in probability we have that
\begin{equation*}
J_{n}(\overline{\beta}_{n},\overline{\mu}_{n})=J_{n}(\beta_{0},\mu_{0}%
)+\sqrt{n}\Lambda_{1\beta}(\mathbf{\beta}_{0},\mu_{0})^{\prime}(\overline
{\mathbf{\beta}}-\mathbf{\beta}_{0})+\sqrt{n}\Lambda_{1\mu}(\mathbf{\beta}%
_{0},\mu_{0})(\overline{\mu}_{n}-\mu_{0})+o_{p}(1). \label{YYY}%
\end{equation*}
\end{lemma}

The proof of  Lemma \ref{tayfa} is based on a small number of
intermediate results, being  Proposition \ref{Hub67-3} the most
important of them. It  may be considered \ \ the U-statistics
version of Lemma 3 of Huber \cite{Huber67}. Since we believe that these
results can be useful in many other situations, we decided to make a
presentation in a general setting.  Consider a sequence
$\mathbf{z}_{i},i\geq1$ of \ i.i.d. random vectors of
dimension $m$ and let $\ \mathbf{\Psi}(\mathbf{z}_{1},\mathbf{z}%
_{2},\mathbf{\theta}):\mathbb{R}^{m}\times\mathbb{R}^{m}\times\Theta
\rightarrow\mathbb{R}^{p}$, where $\Theta\subset\mathbb{R}^{p}.$
Note that here  we are resorting to the same notation already
adopted for the particular case considered above.

 Let
$\mathbf{\Lambda}(\mathbf{\theta})=$\textrm{E}$\Psi(\mathbf{z}_{1}%
,\mathbf{z}_{2,}\mathbf{\theta})$ and assume \ that $\mathbf{\
\Lambda }(\mathbf{\theta}_{0})=0$, for some $\theta_{0}\in R^{p}$. \

Consider
\begin{equation}
\mathbf{Z}_{n}(\mathbf{\theta})=\frac{\left\Vert \sum_{j=1}^{n}\sum_{i=1}%
^{n}[\mathbf{\Psi}(\mathbf{z}_{i},\mathbf{z}_{j},\mathbf{\theta}%
)-\mathbf{\Lambda}(\theta)-\mathbf{\Psi}(\mathbf{z}_{i},\mathbf{z}%
_{j},\mathbf{\theta}_{0})]\right\Vert
}{n^{3/2}+n^{2}||\mathbf{\Lambda
}(\mathbf{\theta})||} \label{Zn}%
\end{equation}
and
\[
U(\mathbf{z}_{1},\mathbf{z}_{2},\mathbf{\theta},d)=\sup_{||\theta^{\ast
}-\mathbf{\theta}||\leq d}\left\Vert \mathbf{\Psi}(\mathbf{z}_{1}%
,\mathbf{z}_{2},\mathbf{\theta}^{\ast})-\mathbf{\Psi}(\mathbf{z}%
_{1},\mathbf{z}_{2},\mathbf{\theta})\right\Vert .
\]
\bigskip We need the following assumptions:

C1. For a fixed $\theta,\mathbf{\Psi}(\mathbf{z}_{1},\mathbf{z}_{2}%
,\mathbf{\theta})$ is measurable and $\mathbf{\Psi}(\mathbf{z}_{1}%
,\mathbf{z}_{2},\mathbf{\theta})$ is separable. For the definition
of separability, see\ Huber \cite{Huber67}.

C2. There exist numbers $b>0,$ $c>0$ and $d_{0}>0$ such that $\ $(i)
$\Lambda(\mathbf{\theta})$ \ is continuously differentiable for
$|\mathbf{|\theta-\theta}_{0}||\leq d_{0}$ and
$\mathbf{\dot{\Lambda}(\theta }_{0})$ is nonsingular, where
$\mathbf{\dot{\Lambda}(\theta)}$ is the
differential matrix of $\Lambda(\mathbf{\theta}),$ (ii) $\mathrm{E}%
U(\mathbf{z}_{1},\mathbf{z}_{2},\mathbf{\theta},d)\leq bd$ if
$||\mathbf{\theta}-\mathbf{\theta}_{0}||+d\leq d_{0}$ and (iii) $\mathrm{E}%
U^{2}(\mathbf{z}_{1},\mathbf{z}_{2},\theta,d)\leq bd$ if
$||\mathbf{\theta }-\mathbf{\theta}_{0}||+d\leq d_{0}.$

C3.
$\mathrm{E}||\mathbf{\Psi}(\mathbf{z}_{1},\mathbf{z}_{2},\mathbf{\theta
}_{0}\mathbf{)||}^{2}<\infty.$

\begin{proposition}
\label{Hub67-3} Suppose \ that \ assumptions C1-C3 hold.\ Then \ we
have
\[
\sup_{||\mathbf{\theta}-\mathbf{\theta}_{0}||\leq d_{0}}\mathbf{Z}%
_{n}(\mathbf{\theta})\rightarrow_{p}0.
\]

\end{proposition}

The proof is similar to that of Lemma 3 in Huber \cite{Huber67}. The only
difference is that all the sums of independent variables need to be
replaced by U-statistics.  Moreover the U-statistics counterparts of
\ $U_{n},V_{n}$ and the right hand side of equation (51) in Huber
\cite{Huber67}, must be approximated by sums of independent random variables
using the Projection Theorem.$\square$

\

Let now $\Psi_{1}(\mathbf{z}_{1},\mathbf{z}_{2},\theta):\mathbb{R}^{m}%
\times\mathbb{R}^{m}\times\Theta\rightarrow\mathbb{R},$ \ where
$\Theta
\subset\mathbb{R}^{p},$ and let $\Lambda_{1}(\mathbf{z}_{1},\mathbf{z}%
_{2},\mathbf{\theta})=\mathrm{E}\Psi_{1}(\mathbf{z}_{1},\mathbf{z}%
_{2},\mathbf{\theta})$. \ Take $\mathbf{\theta}_{0}=(\theta_{01}%
,...,\theta_{0p})$ \ such that $\Lambda_{1}(\mathbf{\theta}_{0})=0$.
Put
\[
U_{1}(\mathbf{z}_{1},\mathbf{z}_{2},\mathbf{\theta},d)=\sup_{||\theta^{\ast
}-\mathbf{\theta}||\leq d}\left\vert \Psi_{1}(\mathbf{z}_{1},\mathbf{z}%
_{2},\theta^{\ast})-\Psi_{1}(\mathbf{z}_{1},\mathbf{z}_{2},\mathbf{\theta
})\right\vert
\]
and%
\[
Z_{n1}(\mathbf{\theta})=\frac{\left\vert
\sum_{j=1}^{n}\sum_{i=1}^{n}[\Psi
_{1}(\mathbf{z}_{i},\mathbf{z}_{j},\mathbf{\theta})-\Lambda_{1}(\mathbf{\theta
})-\Psi_{1}(\mathbf{z}_{i},\mathbf{z}_{j},\mathbf{\theta}_{0})]\right\vert
}{n^{3/2}+n^{2}|\Lambda_{1}(\mathbf{\theta}_{0})|}.
\]
Denote by $\ \ \ \
\mathbf{\dot{\Lambda}}_{1}(\mathbf{\theta})=(\dot{\Lambda
}_{11}(\mathbf{\theta}),...\dot{\Lambda}_{1p}(\mathbf{\theta}))=\partial
\Lambda_{1}(\mathbf{\theta)/\partial}\theta.$

In order to prove a statement
analogous to Proposition \ref{Hub67-3} for the univariate statistics $Z_{n1}(\mathbf{\theta)}$,
the following assumptions will be needed.

D1. For a fixed $\theta$, $\Psi_{1}(\mathbf{z}_{1},\mathbf{z}_{2}%
,\mathbf{\theta})$ is measurable and separable.

D2. There exist numbers $b>0,$ $c>0$ and $d_{0}>0$ such that $\ $(i)
$\Lambda_{1}(\mathbf{\theta})$ is continuously differentiable for
$|\mathbf{|\theta-\theta}_{0}||\leq d_{0}$ and $\dot{\Lambda}_{1}%
(\mathbf{\theta}_{0})$ $\neq0$\ $\mathbf{\ }$ (ii) $\mathrm{E}U_{1}%
(\mathbf{z}_{1},\mathbf{z}_{2},\mathbf{\theta},d)\leq bd$ if
$||\mathbf{\theta
}-\mathbf{\theta}_{0}||+d\leq d_{0}$ and (iii) $\mathrm{E}U_{1}^{2}%
(\mathbf{z}_{1},\mathbf{z}_{2},\mathbf{\theta},d)\leq bd$ if
$||\mathbf{\theta }-\mathbf{\theta}_{0}||+d\leq d_{0}.$

D3 $\mathrm{E}\Psi_{1}^{2}(\mathbf{z}_{1},\mathbf{z}_{2},\mathbf{\theta}%
_{0})<\infty.$

\begin{proposition}
\label{HUB1} Suppose that assumptions \ D1-D3 \ are satisfied$.$ Then%
\begin{equation}
\sup_{||\theta-\mathbf{\theta}_{0}||\leq d_{0}}Z_{n1}(\mathbf{\theta
})\rightarrow_{p}0. \label{Zn1}%
\end{equation}

\end{proposition}

\begin{proof} Let
$\dot{\Lambda}_{1}(\mathbf{\theta})=(\dot{\Lambda}_{11}(\mathbf{\theta
}),...\dot{\Lambda}_{1p}(\mathbf{\theta}))^{\prime}$.
Without loss of generality,  by D2,
we can assume that  $\dot{\Lambda}_{11}(\mathbf{\theta })\neq0$.
For $2\leq i\leq p$, define  $\Psi_{i}(\mathbf{z}_{1},\mathbf{z}_{2},\mathbf{\theta})=\theta_{i}%
-\theta_{0i}$ and consider
$\mathbf{\Psi}(\mathbf{z}_{1},\mathbf{z}_{2},\mathbf{\theta
)=(}\Psi_{1}(\mathbf{z}_{1},\mathbf{z}_{2},\mathbf{\theta),}\Psi
_{2}\mathbf{(\mathbf{z}_{1},\mathbf{z}_{2},\mathbf{\theta),}...},\Psi
_{p}(\mathbf{z}_{1},\mathbf{z}_{2},\mathbf{\theta))}^{\prime}$.
 Doing
$\mathbf{\Lambda(\theta
)=}\mathrm{E}\Psi(\mathbf{z}_{1},\mathbf{z}_{2},\mathbf{\theta),}$
we have
$\mathbf{\Lambda(\theta}_{0}\mathbf{)=0}$ and%

\[
\mathbf{\dot{\Lambda}}(\mathbf{\theta})=\left(
\begin{array}
[c]{cc}%
\dot{\Lambda}_{11}(\mathbf{\theta}) &
\dot{\Lambda}_{12}(\mathbf{\theta
})...\dot{\Lambda}_{1p}(\mathbf{\theta})\\
0 & I_{p-1}%
\end{array}
\right)  .
\]

Then det$(\mathbf{\dot{\Lambda}}(\mathbf{\theta}_{0}))=\dot{\Lambda}%
_{11}(\mathbf{\theta})\neq0$ and it is easy to check that the
remaining assumptions C1-C3 are also satisfied. \ Let
$\mathbf{Z}_{n}(\mathbf{\theta})$ be given by (\ref{Zn}), \ then by
Proposition \ref{Hub67-3}\ we get that
$\sup_{||\mathbf{\theta}-\mathbf{\theta}_{0}||\leq d_{0}}\mathbf{Z}%
_{n}(\mathbf{\theta})\rightarrow_{p}0.$ This implies (\ref{Zn1}).
\end{proof}

\begin{proposition}
\label{corlem} Suppose the same assumptions as in Proposition
\ref{HUB1} \ and let $\overline{\mathbf{\theta}}_{n}$ be a sequence
of estimates of
$\ \mathbf{\theta}_{0}$ such that $n^{1/2}|||\overline{\mathbf{\theta}}%
_{n}-\mathbf{\theta}_{0}||=O_{p}(1).$ Then%
\begin{equation}
\frac{1}{n^{3/2}}\sum_{j=1}^{n}\sum_{i=1}^{n}\Psi_{1}(\mathbf{z}%
_{i},\mathbf{z}_{j},\overline{\mathbf{\theta}}_{n})=\frac{1}{n^{3/2}}%
\sum_{j=1}^{n}\sum_{i=1}^{n}\Psi_{1}(\mathbf{z}_{i},\mathbf{z}_{j}%
,\mathbf{\theta}_{0})+\dot{\Lambda}_{1}(\mathbf{\theta}_{0})^{\prime}%
(n^{1/2}\
(\overline{\mathbf{\theta}}_{n}-\mathbf{\theta}_{0}))+o_{p}(1).
\label{lecorlem}%
\end{equation}
\end{proposition}

\begin{proof}
By Proposition \ref{HUB1} \ we have%
\begin{equation}
Z_{n1}(\overline{\mathbf{\theta}}_{n})=\frac{\left\vert \sum_{j=1}^{n}%
\sum_{i=1}^{n}[\Psi_{1}(\mathbf{z}_{i},\mathbf{z}_{j},\overline{\mathbf{\theta
}}_{n})-\Lambda_{1}(\overline{\mathbf{\theta}}_{n})-\Psi_{1}(\mathbf{z}%
_{i},\mathbf{z}_{j},\mathbf{\theta}_{0})]\right\vert }{n^{3/2}+n^{2}%
|\Lambda_{1}(\overline{\mathbf{\theta}}_{n})|}\rightarrow_{p}0. \label{con0}%
\end{equation}
Using the Mean Value Theorem\ \ we get%
\begin{align*}
&  \frac{\left\vert \sum_{j=1}^{n}\sum_{i=1}^{n}[\Psi_{1}(\mathbf{z}%
_{i},\mathbf{z}_{j},\overline{\mathbf{\theta}}_{n})-\dot{\Lambda}%
_{1}(\mathbf{\theta}_{n}^{\ast})^{\prime}(\overline{\mathbf{\theta}}%
_{n}-\mathbf{\theta}_{0})-\Psi_{1}(\mathbf{z}_{i},\mathbf{z}_{j}%
,\mathbf{\theta}_{0})]\right\vert }{n^{3/2}+n^{2}\left\vert \dot{\Lambda}%
_{1}(\mathbf{\theta}_{n}^{\ast})^{\prime}(\overline{\mathbf{\theta}}%
_{n}-\mathbf{\theta}_{0})\right\vert }\\
&  =\frac{\left\vert \sum_{j=1}^{n}\sum_{i=1}^{n}[\Psi_{1}(\mathbf{z}%
_{i},\mathbf{z}_{j},\overline{\mathbf{\theta}}_{n})-\dot{\Lambda}%
_{1}(\mathbf{\theta}_{n}^{\ast})^{\prime}(\overline{\mathbf{\theta}}%
_{n}-\mathbf{\theta}_{0})-\Psi_{1}(\mathbf{z}_{i},\mathbf{z}_{j}%
,\mathbf{\theta}_{0})]\right\vert }{n^{3/2}(1+\left\vert \dot{\Lambda}%
_{1}(\mathbf{\theta}_{n}^{\ast})^{\prime}n^{1/2}(\overline{\mathbf{\theta}%
}_{n}-\mathbf{\theta}_{0})\right\vert )},
\end{align*}
where
$\mathbf{\theta}_{n}^{\ast}\rightarrow_{p}\mathbf{\theta}_{0}.$
Since
$\dot{\Lambda}(\mathbf{\theta}_{n}^{\ast})n^{1/2}(\overline{\mathbf{\theta}%
}_{n}-\mathbf{\theta}_{0})$ is bounded in probability, \ (\ref{con0}) implies%
\[
\frac{1}{n^{3/2}}\sum_{j=1}^{n}\sum_{i=1}^{n}[\Psi_{1}(\mathbf{z}%
_{i},\mathbf{z}_{j},\overline{\mathbf{\theta}}_{n})-\dot{\Lambda}%
_{1}(\mathbf{\theta}_{n}^{\ast})^{\prime}(\overline{\mathbf{\theta}}%
_{n}-\mathbf{\theta}_{0})-\Psi_{1}(\mathbf{z}_{i},\mathbf{z}_{j}%
,\mathbf{\theta}_{0})]\ \rightarrow_{p}0,
\]
and so
\[
\frac{1}{n^{3/2}}\sum_{j=1}^{n}\sum_{i=1}^{n}\Psi_{1}(\mathbf{z}%
_{i},\mathbf{z}_{j},\overline{\mathbf{\theta}}_{n})=\frac{1}{n^{3/2}}%
\sum_{j=1}^{n}\sum_{i=1}^{n}\Psi_{1}(\mathbf{z}_{i},\mathbf{z}_{j}%
,\mathbf{\theta}_{0})+\dot{\Lambda}_{1}(\mathbf{\theta}_{n}^{\ast}%
)(n^{1/2}\
(\overline{\mathbf{\theta}}_{n}-\mathbf{\theta}_{0}))+o_{p}(1).
\]
Finally, using the continuity of \ $\dot{\Lambda}$ at
$\mathbf{\theta}_{0}$, the order of convergence of
$\overline{\theta}_{n}$, and the fact that
$\mathbf{\theta}_{n}^{\ast}\rightarrow_{p}\mathbf{\theta}_{0}$, \ we
get (\ref{lecorlem}). \end{proof}

In the following Proposition we give closed formulas \ for
$\Lambda_{1\mathbf{\beta }}(\mathbf{\beta}_{0},\mu_{0}\mathbf{)}$
and $\Lambda_{1\mu}(\mathbf{\beta }_{0},\mu_{0})$, which are part of the expansion stated in Lemma \ref{tayfa}.

\begin{proposition}
\label{caldev} We have%
\[
\Lambda_{1\mathbf{\beta}}(\mathbf{\beta}_{0},\mu_{0})=2\mathrm{E}[a_{1}%
k_{0}(-g(\mathbf{x}_{2},\mathbf{\beta}_{0})+\mu_{0})({g}(\mathbf{x}%
_{2},\mathbf{\beta}_{0})-{g}(\mathbf{x}_{1},\mathbf{\beta}_{0}))]
\]
and%
\[
\Lambda_{1\mu}(\mathbf{\beta}_{0},\mu_{0}\mathbf{)=-}2\eta
f_{0}(\mu_{0}).
\]

\end{proposition}

\begin{proof} Let $\Delta($\textbf{$x$}$_{i},\mathbf{\beta})=g(\mathbf{x}%
_{i},\mathbf{\beta})-g(\mathbf{x}_{i},\mathbf{\beta}_{0}).$ Then
\begin{align*}
\Lambda_{1}(\mathbf{\beta,\mu})  &  =\eta\mathrm{E}(\text{sign}(u_{1}%
+g(\mathbf{x}_{2},\mathbf{\beta}_{0})+\Delta(\mathbf{x}_{2},\mathbf{\beta
})-\Delta(\mathbf{x}_{1},\mathbf{\beta})-\mu|a_{1}=1)\\
&  =\eta\mathrm{E}(\mathrm{E}(\text{sign}(u_{1}+g(\mathbf{x}_{2}%
,\mathbf{\beta}_{0})+\Delta(\mathbf{x}_{2},\mathbf{\beta})-\Delta
(\mathbf{x}_{1},\mathbf{\beta})-\mu|a_{1}=1,\mathbf{x}_{1},\mathbf{x}%
_{2}\mathbf{))}\\
&  =\eta\mathrm{E}((1-2K_{0}(-g(\mathbf{x}_{2},\mathbf{\beta}_{0}%
)-\Delta(\mathbf{x}_{2},\mathbf{\beta})+\Delta(\mathbf{x}_{1},\mathbf{\beta
})+\mu)|a_{1}=1).
\end{align*}
Differentiating the last equation we get%
\begin{align*}
{\Lambda_{1}}_{\mathbf{\beta}}(\mathbf{\beta},\mu)  &  =-2\eta\mathrm{E}%
[k_{0}(-g(\mathbf{x}_{2},\mathbf{\beta}_{0})-\Delta(\mathbf{x}_{2}%
,\mathbf{\beta})+\Delta(\mathbf{x}_{1},\mathbf{\beta})+\mu)(-{g}%
(\mathbf{x}_{2},\mathbf{\beta})+{g}(\mathbf{x}_{1},\mathbf{\beta
}))]|a=1)\nonumber\\
&
=-2\eta\mathrm{E}[k_{0}(-g(\mathbf{x}_{2},\mathbf{\beta}_{0})-\Delta
(\mathbf{x}_{2},\mathbf{\beta})+\Delta(\mathbf{x}_{1},\mathbf{\beta}%
)+\mu)|a=1)\nonumber\\
&  =-2\eta\mathrm{E}[k_{0}(-g(\mathbf{x}_{2},\mathbf{\beta}_{0})+\mu_{0}%
)(-{g}(\mathbf{x}_{2},\mathbf{\beta}_{0})+{g}(\mathbf{x}_{1},\mathbf{\beta
}_{0}))]|a_{1}=1)\nonumber\\
&  =2\mathrm{E}[a_{1}k_{0}(-g(\mathbf{x}_{2},\mathbf{\beta}_{0})+\mu_{0}%
)({g}(\mathbf{x}_{2},\mathbf{\beta}_{0})-{g}(\mathbf{x}_{1},\mathbf{\beta}%
_{0}))]
\end{align*}
and
\begin{align*}
{\Lambda_{1}}_{\mathbf{\mu}}(\mathbf{\beta}_{0},\mu_{0})  &  =-2\eta
\mathrm{E}[k_{0}(-g(\mathbf{x}_{2},\mathbf{\beta}_{0})+\mu_{0})]\nonumber\\
&  =-2\eta f_{0}(\mu_{0}).
\end{align*}

These prove the Proposition. \end{proof}

\textsc{Proof of Lemma \ref{tayfa}:} By Proposition \ref{corlem}, we \ only
need to verify that under the assumptions of Theorem 3, considering
$\mathbf{\theta
=}(\mathbf{\beta,}\mu)$,%

\[
\Psi_{1}(\mathbf{z}_{1},\mathbf{z}_{2},\mathbf{\theta})=a_{1} \text{sign}%
(y_{1}-g(\mathbf{x}_{1},\mathbf{\beta})+g(\mathbf{x}_{2},\mathbf{\beta})-\mu).
\]
and
\[
U_{1}(\mathbf{z}_{1},\mathbf{z}_{2},\mathbf{\theta},d)=\sup_{||\theta^{\ast
}-\theta||\leq d}|\Psi_{1}(\mathbf{z}_{1},\mathbf{z}_{2},\mathbf{\theta}%
^{\ast})-\Psi_{1}(\mathbf{z}_{1},\mathbf{z}_{2},\mathbf{\theta})|,
\]
then, assumptions D1-D3 are satisfied.

Assumptions D1 and D3 follow immediately. Assumption D2(i) follows
from Proposition  \ref{caldev} and the fact that $f_{0}(\mu_{0})>0$.

We now prove D2 (ii) and (iii). Take $d_{0}=\delta$ as in A0, then if we put%
\begin{equation}
\omega(\mathbf{x})=\sup_{\left\Vert \mathbf{\beta}-\mathbf{\beta}%
_{0}\right\Vert \leq\delta}\left\Vert
{\dot{g}}(\mathbf{x},\mathbf{\beta })\right\Vert ,
\end{equation}
we have $\mathrm{E}\omega(\mathbf{x}_{1})<\infty$. To prove D2 (i)
and (ii), we have to show that there exist $K_{1}$ and $K_{2}$ such
that for all $\mathbf{\theta}$ \ and $d$ with
$||\theta-\theta_{0}||+d\leq d_{0}$, we have
\begin{equation}
\mathrm{E}U_{1}^{i}(\mathbf{z}_{1},\mathbf{z}_{2},\mathbf{\theta},d)\leq
K_{i}d\;,\text{ }i=1,2.\label{doudes}%
\end{equation}
\ For that purpose, we can write
\begin{align}
U_{1}(\mathbf{z}_{1},\mathbf{z}_{2},\mathbf{\theta,}d) &  \leq\sup
_{||\mathbf{\beta}^{\ast}-\mathbf{\beta}||\leq
d,|\mu^{\ast}-\mu|\leq
d}|\text{sign}(u_{1}+g(\mathbf{x}_{1},\mathbf{\beta}_{0})-g(\mathbf{x}%
_{1},\mathbf{\beta}^{\ast})+g(\mathbf{x}_{2},\mathbf{\beta}^{\ast})-\mu^{\ast
}) \nonumber\\
&  -\text{sign}(u_{1}+g(\mathbf{x}_{1},\mathbf{\beta}_{0})-g(\mathbf{x}%
_{1},\mathbf{\beta})+g(\mathbf{x}_{2},\mathbf{\beta})-\mu)|.\label{HH1}%
\end{align}

Then if $||\theta-\theta_{0}||+d\leq d_{0},$ $\ $and
$||\mathbf{\theta}^{\ast }-\mathbf{\theta}||\leq d$ \ we get that
$||\mathbf{\theta}^{\ast }-\mathbf{\theta}_{0}||\leq d_{0}$, and
therefore $||\mathbf{\beta}^{\ast }-\mathbf{\beta}_{0}||\leq d_{0}$
too. Note also  that
\begin{align}
&  \left\vert
(g(\mathbf{x}_{1},\beta_{0})-g(\mathbf{x}_{1},\mathbf{\beta
}^{\ast})+g(\mathbf{x}_{2},\mathbf{\beta}^{\ast})-\mu^{\ast})-(g(\mathbf{x}%
_{1},\beta_{0})-g(\mathbf{x}_{1},\mathbf{\beta})+g(\mathbf{x}_{2}%
,\mathbf{\beta})-\mu)\right\vert \nonumber\\
&
\leq(\omega(\mathbf{x}_{1})+\omega(\mathbf{x}_{2}))||\mathbf{\beta}^{\ast
}-\mathbf{\beta}||+|\mu^{\ast}-\mu|\leq(\omega(\mathbf{x}_{1})+\omega
(\mathbf{x}_{2})+1)d. \label{HH2}%
\end{align}

Let $z=$ $(\omega(\mathbf{x}_{1})+\omega(\mathbf{x}_{2})+1)$,
$v=|u|$ and $w=v/z.$ \ The left hand side of (\ref{HH1}) is
different from 0 when the two arguments of the sign function have
different signs. By (\ref{HH2}) \ this occurs \ only if
$|u_{1}|\leq(\omega(\mathbf{x}_{1})+\omega (\mathbf{x}_{2})+1)d.$
Then we can write
\[
U_{1}(\mathbf{z}_{1},\mathbf{z}_{2},\mathbf{\beta},\mu\mathbf{,}d)\leq2I(w\leq
d).
\]
Observe that   $\mathrm{E}z<\infty$ and \ that the density of $v$ is
given by $f_{v}(v)=k_{0}(v)+k_{0}(-v)$, which is bounded by $2\sup
k_{0}$. Then, since
\ the density of $\ w$ is%
\[
f_{w}(w)=%
{\displaystyle\int\limits_{0}^{\infty}}
zf_{v}(wz)dF_{z}\leq2\sup k_{0}%
{\displaystyle\int\limits_{0}^{\infty}}
zdF_{z}=2\sup k_{0}\mathrm{E}(z)
\]
$\ $ we get
\[
\mathrm{E}U_{1}(\mathbf{z}_{1},\mathbf{z}_{2},\mathbf{\beta},\mu
\mathbf{,}d)\leq2P(w\leq d)\leq4\sup k_{0}\mathrm{E}(z)d
\]
and%
\[
\mathrm{E}U_{1}^{2}(\mathbf{z}_{1},\mathbf{z}_{2},\mathbf{\beta},\mu
\mathbf{,}d)\leq4P(w\leq d)\leq8\sup k_{0}\mathrm{E}(z)d,
\]
and so (\ref{doudes}) holds with $K_{1}=4\sup k_{0}$\textrm{E}$(z)$
and $K_{2}=8\sup k_{0}$\textrm{E}$(z).$

The expansion  obtained in Lemma \ref{tayfa} requires that $\sqrt
{n}(\overline{\mathbf{\theta}}_{n}-\mathbf{\theta}_{0})=O_{P}(1)$. The
following Lemma shows that \ $\widehat{\mathbf{\theta}}_{n}=(\widehat
{\mathbf{\beta}}_{n},\widehat{\mu}_{n})$ satisfies this condition.

\begin{lemma}
\label{jacot} Under the assumptions of Theorem  \ \ref{median} we
have that
\newline(a) $n^{1/2}(\widehat{\mu}_{n}-\mu_{0})$ is bounded in probability,
\newline(b) $J_{n}(\widehat{\mathbf{\beta}}_{n},\widehat{\mu}_{n}%
)\rightarrow^{p}0.$
\end{lemma}

\begin{proof} Let
\[
J_{n}^{\ast}(\mathbf{\beta,}\mu)=\frac{1}{n^{3/2}}\sum_{j=1}^{n}\sum_{i=1}%
^{n}\Psi(\mathbf{z}_{i},\mathbf{z}_{j},\mathbf{\beta}_{0},\mu_{0}%
))+n^{1/2}\mathbf{\Lambda_{1}}_{\mathbf{\beta}}(\mathbf{\beta}_{0},\mu
_{0})^{\prime}(\mathbf{\beta}-\mathbf{\beta}_{0})+n^{1/2}\Lambda
_{1\mathbf{\mu}}(\mathbf{\beta}_{0},\mu_{0})(\mu-\mu_{0}).
\]
Take $\varepsilon>0$ and let $\widetilde{\mu}_{1n}$ $\ $\ and $\widetilde{\mu
}_{2n}$ be defined by
\[
J_{n}^{\ast}(\widehat{\mathbf{\beta}}_{n}\mathbf{,}\widetilde{\mu}%
_{1n})=\varepsilon\text{ and }J_{n}^{\ast}(\widehat{\mathbf{\beta}}%
_{n}\mathbf{,}\widetilde{\mu}_{2n})=-\varepsilon.\]

By Lemma \ref{caldev} \ $\Lambda_{1\mu}(\mathbf{\beta}_{0},\mu
_{0})\neq0,$ and it holds that
\[
n^{1/2}(\widetilde{\mu}_{1n}-\mu_{0})=-\frac{1}{\Lambda_{1\mu}(\mathbf{\beta
}_{0},\mu_{0})}\left[  \frac{1}{n^{3/2}}\sum_{j=1}^{n}\sum_{i=1}^{n}%
\Psi(\mathbf{z}_{i},\mathbf{z}_{j},\mathbf{\beta}_{0},\mu_{0})+n^{1/2}%
\mathbf{\Lambda_{1}}_{\mathbf{\beta}}(\mathbf{\beta}_{0},\mu_{0})^{\prime
}(\widehat{\mathbf{\beta}}_{n}-\mathbf{\beta}_{0})-\varepsilon\right]  .
\]
\ Both the first and second terms on the right hand side are bounded
in probability, the former by the Central Limit Theorem for
U-statistics and the later by Assumption A2 and the Central Limit
Theorem. Then $n^{1/2}(\widetilde{\mu}_{1n}-\mu_{0})$ is
bounded in probability. \ Thereafter, by Lemma \ref{tayfa} we get%

\[
J_{n}(\widehat{\mathbf{\beta}}_{n},\widetilde{\mu}_{1n})=\varepsilon
+o_{p}(1).
\]
Similarly we can prove that $n^{1/2}(\widetilde{\mu}_{2n}-\mu_{0})$
is bounded in probability and that
\[
J_{n}(\widehat{\mathbf{\beta}}_{n},\widetilde{\mu}_{2n})=-\varepsilon
+o_{p}(1).
\]
Then since $J(\beta,\mu)$ is nonincreasing in $\mu,$ \ by a property of the
median we get that
\[
\lim_{n\rightarrow\infty}\mathrm{P}(\widetilde{\mu}_{1n}\leq\widehat{\mu}%
_{n}\leq\widetilde{\mu}_{2n})=1,
\]
and therefore $n^{1/2}(\widehat{\mu}_{n}-\mu_{0})$ is bounded in probability.
We also have that
\[
\mathrm{P}(J_{n}(\widehat{\mathbf{\beta}}_{n},\widetilde{\mu}_{2n})\leq
J_{n}(\widehat{\mathbf{\beta}}_{n},\widehat{\mu}_{n})\leq J_{n}(\widehat
{\mathbf{\beta}}_{n},\widetilde{\mu}_{1n}))\rightarrow1
\]
and therefore
\[
\mathrm{P}(-2\varepsilon\leq J_{n}(\widehat{\mathbf{\beta}}_{n},\widehat{\mu
}_{n})\leq2\varepsilon)\rightarrow1.
\]
Since this holds for all $\varepsilon>0,$ part (b) of the Lemma \ is proved.
\end{proof}

\textsc{Proof of Theorem \ref{median}.}

(a) To prove this part of the Theorem it suffices to show that $T_{\text{med}%
}$ is \ weakly continuous at $F_{0}.$ \ Take $\varepsilon>0$ and $y_{1\text{
}}$ , $y_{2}$ continuity points of $\ F_{0}$ $\ $\ such that \ $\mu
_{0}-\varepsilon<y_{1}<\mu_{0}<y_{2}<\mu_{0}+\varepsilon$. Since $F_{0}$ is
continuous and strictly increasing at $F_{0}$ \ we have that $F_{0}(\mu
_{0})=0.5$ and there exists $\delta>0$ such that\ $F_{0}(y_{1})<0.5-\delta
<0.5+\delta<F_{0}(\mu_{2}).$ \ Suppose that $F_{n}\rightarrow_{w}F_{0},$ then
there exists $\ n_{0}$ such that for $n\geq n_{0}$ we have $F_{n}%
(y_{1})<0.5-\delta$ and $F_{n}(y_{2})>0.5+.\delta.$ This proves that $n\geq
n_{0}$ implies \ that $\ \mu_{0}-\varepsilon<y_{1}\leq T_{\text{med}}%
(F_{n})\leq y_{2}<\mu_{0}+\varepsilon.$\

(b) \ Since $n^{1/2}(\widehat{\mu}_{n}-\mu_{0})$ and $n^{1/2}(\widehat
{\mathbf{\beta}}_{n}-\mathbf{\beta}_{0})$\ are bounded in probability, \ by
Lemma\ \ref{tayfa} \ \ we get
\[
J_{n}(\widehat{\mathbf{\beta}}_{n},\widehat{\mu}_{n})=J_{n}(\mathbf{\beta}%
_{0},\mu_{0})+\sqrt{n}\Lambda_{1\beta}(\mathbf{\beta}_{0},\mu_{0}%
)(\widehat{\mathbf{\beta}}_{n}-\mathbf{\beta}_{0})+\sqrt{n}\Lambda_{1\mu
}(\beta_{0},\mu_{0})(\widehat{\mu}_{n}-\mu_{0})\}+o_{p}(1),
\]
and using Lemma \ref{jacot} (b) we get
\begin{align*}
\sqrt{n}\{\widehat{\mu}_{n}-\mu_{0}\}  &  =-\frac{1}{\Lambda_{1\mu}(\beta
_{0},\mu_{0})\ }\left\{  \Lambda_{1\mathbf{\beta}}(\mathbf{\beta}_{0},\mu
_{0})^{\prime}n^{1/2}(\widehat{\mathbf{\beta}}_{n}-\mathbf{\beta}_{0}%
)+J_{n}(\mathbf{\beta}_{0},\mu_{0})\right\}  +o_{p}(1)\\
&  =d_{n}^{\ast}+n^{1/2}\mathbf{c}_{n}^{\ast\prime}(\widehat{\mathbf{\beta}%
}_{n}-\mathbf{\beta}_{0})+o_{p}(1),
\end{align*}
where
\[
\mathbf{c}_{n}^{\ast}=\frac{1}{_{n^{2}}}\sum_{j=1}^{n}\sum_{i=1}^{n}%
a_{i}\,I_{L}^{\prime}(y_{i}-g({\mathbf{\beta}},\mathbf{x}_{i})+{g}%
({\mathbf{\beta}},\mathbf{x}_{j}))\,\{\dot{g}({\mathbf{\beta}},\mathbf{x}%
_{j})-\dot{g}({\mathbf{\beta}},\mathbf{x}_{i})\}
\]
and%
\[
d_{n}^{\ast}=\frac{1}{n^{3/2}}%
{\displaystyle\sum_{i=1}^{n}}
{\displaystyle\sum_{j=1}^{n}}
a_{i}I_{T_{med}}(u_{i}+g(\mathbf{\beta}_{0},\mathbf{x}_{j})).
\]
\ \bigskip Then it suffices to show that%

\[
V_{n}^{\ast}=d_{n}^{\ast}+\mathbf{c}_{n}^{\ast\prime}n^{1/2}(\widehat
{\mathbf{\beta}}_{n}-\mathbf{\beta}_{0})\rightarrow N(0,\tau^{2}).
\]
The proof of this result is similar to that \ of $d_{n}+\mathbf{c}_{n}%
^{\prime}n^{1/2}(\widehat{\mathbf{\beta}}_{n}-\mathbf{\beta}_{0}%
)\rightarrow^{d}$N$(0,\tau^{2})$ in Theorem \ref{asym}. $\square$\bigskip\

\textsc{Proof of Theorem \ref{BDP}}. Let\textbf{\ }$\mathbf{W}$ be as \ in
(\ref{Wcom}). We have to show that given $t<n\varepsilon_{3}$ and
$s<m\varepsilon_{3},$ \ there exists $K$ such that \ for any sample
$\mathbf{W}^{\ast}\in\mathcal{W}_{ts},$\ we have that $|T_{L}(\widehat{F}%
_{n}^{\ast})|\leq K,$ where $\ \widehat{F}_{n}^{\ast}$ is the distribution
constructed as in (\ref{Gnsom}), based on $\mathbf{W}^{\ast}$. According to
the definition of $\varepsilon_{U}^{\ast}(T_{L})$, it suffices to show that
there exists $M$ such that for any $\mathbf{W}^{\ast}\in\mathcal{W}_{ts}$ we
have that the corresponding $\widehat{F}_{n}^{\ast}$ satisfies
\begin{equation}
\mathrm{P}_{\widehat{F}_{n}^{\ast}}(|y|\leq M)>1-\varepsilon_{2}. \label{Ihin}%
\end{equation}

Let
\[
\mathcal{Z}_{s}=\{\mathbf{Z}^{\ast}=\{(\mathbf{x}_{i}^{\ast},y_{i}^{\ast
}):i\in A\}:%
{\displaystyle\sum_{i\in A}}
I\{(\mathbf{x}_{i}^{\ast},y_{i}^{\ast})\neq(\mathbf{x}_{i},y_{i})\}\leq s\}.
\]
Since $s/m<\varepsilon_{1}$ we can find $M_{1}$ such that
\begin{equation}
\sup_{\mathbf{Z}^{\ast}\in\mathcal{Z}_{s}}||\widetilde{\mathbf{\beta}}%
_{m}(\mathbf{Z}^{\ast})||\leq M_{1}, \label{cot0}%
\end{equation}
and then we can find $M$ such that
\begin{equation}
\sup_{1\leq j\leq n}\sup_{|\mathbf{|\beta}||\leq M_{1}}|g(\mathbf{x}%
_{j},\mathbf{\beta}|\leq M/2 \label{cot1}%
\end{equation}
and%
\begin{equation}
\sup_{i\in A}\sup_{||\beta||\leq M_{1}}|y_{i}-g(\mathbf{x}_{i},\mathbf{\beta
})|\leq M/2. \label{cot2}%
\end{equation}

Given $\ \mathbf{W}^{\ast}\in\mathcal{W}_{t,s}$, if $\widehat{\beta}_{n}%
^{\ast}=\widetilde{\beta}_{m}(\mathbf{Z}^{\ast})$, with $\mathbf{Z}^{\ast}%
\in\mathcal{Z}_{s}$. Consider $B=\{j:$ $1\leq j\leq n,\mathbf{x}%
_{j}=\mathbf{x}_{j}^{\ast}\}$ and $C=\{i\in A:$ $(\mathbf{x}_{j}%
,y_{j})=(\mathbf{x}_{j}^{\ast},y_{j}^{\ast})\}$. Then $\#B>(1-\varepsilon
_{3})n$ and $\#C>(1-\varepsilon_{3})m$. For $1\leq j\leq n,$ $i\in A,$ put
$\widehat{y}_{ij}^{\ast}=g(\mathbf{x}_{j}^{\ast},\widehat{\mathbf{\beta}}%
_{n}^{\ast})+(y_{i}^{\ast}-g(\mathbf{x}_{i}^{\ast},\widehat{\mathbf{\beta}%
}_{n}^{\ast})).$ Then, when $\ j\in B$ and $\ i\in C$, by (\ref{cot0}),
(\ref{cot1}) and (\ref{cot2}), we have that $|\widehat{y}_{ij}^{\ast}|\leq M$
and so
\[
\#\{(i,j):|\widehat{y}_{ij}^{\ast}|\leq M\}>mn(1-\varepsilon_{3})^{2}%
\geq(1-\varepsilon_{2})mn.
\]
Since there are $mn$ pairs $(i,j)$ \ subindexing \ \ \ $\widehat{y}_{ij}%
^{\ast}$, we get that \ $\mathrm{P}_{\widehat{F}_{n}^{\ast}}(|y|$ $\leq
M)>1-\varepsilon_{2}$ and then (\ref{Ihin}) holds. $\square$\bigskip

\textsc{Proof of Theorem \ref{asslreg}}. The proof of this Theorem
is essentially based on Theorem 7 \ of Fasano et. at. \cite{Fasano
et al}. As is mentioned in Section \ref{results}, \ if
$(\mathbf{x}_{i},y_{i})$ has distribution
$G_{0}^{\ast},$ then (\ref{modelogral}) is satisfied with $\mathbf{x}%
_{i}^{\ast}$ having distribution $Q_{0}^{\ast}$ and $u_{i}^{\ast}$
with distribution $K_{0}.$ Moreover, since by Lemma \ref{g-c}
$G_{n}^{\ast }\rightarrow_{w}G_{0}^{\ast}$ , by parts (i), (ii) and
(iii) of Theorem 7 of \cite{Fasano et al} with $G_{0}$ replaced by
$G_{0}^{\ast}$, we get part (i) of the present Theorem.

We now prove (ii). We start proving that for any function $d$ such the
\textrm{$E$}$_{G_{0}^{\ast}}\left\vert d(\mathbf{x},y)\right\vert <\infty$, we
have that
\begin{equation}
\mathrm{E}_{G_{n}^{\ast}}d(\mathbf{x},y)\rightarrow\mathrm{E}_{G_{0}^{\ast}%
}d(\mathbf{x},y)\,\text{\ a.s.} \label{con3}%
\end{equation}
\ Since%
\[
\mathrm{E}_{G_{n}^{\ast}}d(\mathbf{x},y)=\frac{%
{\displaystyle\sum\limits_{i=1}^{n}}
a_{i}d(\mathbf{x}_{i},y_{i})}{%
{\displaystyle\sum\limits_{i=1}^{n}}
a_{i}}\ =\frac{1}{\eta_{n}}\frac{1}{n}%
{\displaystyle\sum\limits_{i=1}^{n}}
a_{i}d(\mathbf{x}_{i},y_{i})
\]

and $\eta_{n}\rightarrow\eta,$ by the Law of Large Numbers \ we have that
\textrm{$E$}$_{G_{n}^{\ast}}d(\mathbf{x},y)\rightarrow$\textrm{$E$}%
$a_{1}d(\mathbf{x}_{1},y_{1})/\eta$ a.s. Since $\mathrm{E}a_{1}d(\mathbf{x}%
_{1},y_{1})/\eta$ $=$\textrm{$E$}$_{G_{0}^{\ast}}d(\mathbf{x},y)$ , we obtain
(\ref{con3}).

Put now $\mathbf{T=(T}_{S},\mathbf{T}_{MM},S)$ and let $I_{\mathbf{T,}%
G_{0}^{\ast}}(\mathbf{x,}y)$ be its influence function at $G_{0}^{\ast}.$ We
now prove that
\[
n^{1/2}\mathrm{E}_{G_{n}^{\ast}}I_{\mathbf{T,}G_{0}^{\ast}}(\mathbf{x,}%
y)\rightarrow_{d}H,
\]
where $H$ is a multivariate normal distribution. This follows by applying the
Central Limit Theorem from%
\[
n^{1/2}\mathrm{E}_{G_{n}^{\ast}}I_{\mathbf{T,}G_{0}^{\ast}}(\mathbf{x,}%
y)=\frac{1}{\eta_{n}}\frac{1}{n^{1/2}}%
{\displaystyle\sum\limits_{i=1}^{n}}
a_{i}I_{\mathbf{T,}G_{0}^{\ast}}(\mathbf{x}_{i},y_{i}),
\]
the facts that \textrm{$E$}$_{G_{0}^{\ast}}I_{\mathbf{T,}G_{0}^{\ast}%
}(\mathbf{x},y)=0,$ and the fact\ that under $G_{0}^{\ast},$ the influence
function $I_{\mathbf{T,}G_{0}^{\ast}}(\mathbf{x},y)$ has\ finite second
moments. Then all the conditions required to apply parts (iv) and (v) of
Theorem 7 of Fasano et al. \cite{Fasano et al} are satisfied. Then%
\[
n^{1/2}(\mathbf{T_{MM,\beta}}(G_{n}^{\ast})-\mathbf{\beta}_{0})=n^{-1/2}%
\sum_{i=1}^{n}a_{i}\frac{1}{\mathrm{E}[a_{1}]}I_{\mathbf{T_{MM,\mathbf{\beta}%
}},G_{0}^{\ast}}(\mathbf{x}_{i},y_{i})+o_{P}(1).
\]
Finally, using \ the expression for $I_{\mathbf{T_{MM,\mathbf{\beta}}}%
,G_{0}^{\ast}}$ derived in Fasano et. al. \cite{Fasano et al}, we obtain part (ii) of the
Theorem. Part (iii) is an immediate consequence of the fact \ that in this
case $e_{01}=0.\square$\ \bigskip

\textsc{Proof of Theorem \ref{assloc}.} Part (i) \ follows from parts (i),
(ii) and (iii) of Theorem 8 of Fasano et al.\cite{Fasano et al}. Let $\mathbf{T}^{L}(F)$
be the complete functional
\[
\mathbf{T}^{L}\left(  F\right)  =\left(  {T}_{{S}}^{L}\left(  F\right)
,{T}_{{MM}}^{L}\left(  F\right)  ,S_{L}\left(  F\right)  \right)  .
\]
\ Since $\{\widehat{F}_{n}\}$ is a sequence of random distribution with finite
support converging a.s. to $F_{0}$, \ by part (iv) \ of Theorem 8 of Fasano et
al. \cite{Fasano et al} we get that $\mathbf{T}^{L}$ is weakly differentiable at
$\{\widehat{F}_{n}\}$ a.s., and so
\begin{equation}
\mathbf{T}^{L}(\widehat{F}_{n})-\mathbf{T}^{L}(F_{0})=\mathrm{E}_{\widehat
{F}_{n}}I_{\mathbf{T}^{L},F_{0}}(y\mathbf{)+}o\left(  \left\Vert
\mathrm{E}_{\widehat{F}_{n}}I_{\mathbf{T}^{L},F_{0}}(y\mathbf{)}\right\Vert
\right)  , \label{expb}%
\end{equation}
where $I_{\mathbf{T}^{L},F_{0}}$ is the influence function of $\mathbf{T}^{L}$
at $F_{0}$.

We prove now that $n^{1/2}\mathrm{E}_{\widehat{F}_{n}}I_{\mathbf{T}^{L},F_{0}%
}(y\mathbf{)}$ is bounded in probability. Using a Taylor expansion, we get%
\begin{equation}
\sqrt{n}\,\mathrm{E}_{\widehat{F}_{n}}I_{T^{L},F_{0}}(y)=\frac{1}{\eta_{n}%
}\left\{  D_{n}+\mathbf{C}_{n}^{\text{'}}n^{1/2}(\widehat{\mathbf{\beta}}%
_{n}-\mathbf{\beta}_{0})\right\}  , \label{TE1}%
\end{equation}
where%
\begin{equation}
D_{n}={n^{-3/2}}%
{\displaystyle\sum_{i=1}^{n}}
{\displaystyle\sum_{j=1}^{n}}
a_{i}I_{T^{L},F_{0}}(u_{i}+g(x_{j},\mathbf{\beta}_{0})), \label{TE2}%
\end{equation}

\begin{equation}
\mathbf{C}_{n}=\frac{1}{n^{2}}\sum_{i=1}^{n}\sum_{j=1}^{n}h(\mathbf{x}%
_{i},y_{i},a_{i},\mathbf{x}_{j},\beta_{n}^{\ast}),\ \label{TE3}%
\end{equation}
$\mathbf{\beta}_{n}^{\ast}$ between $\widehat{\mathbf{\beta}}_{n}$ and
$\mathbf{\beta}_{0}$ and
\[
h(\mathbf{x}_{i},y_{i},a_{i},\mathbf{x}_{j},\beta)=a_{i}\,I_{{T^{L},F_{0}}%
}^{\prime}(y_{i}-g(x_{i},{\mathbf{\beta}})+g(x_{j},{\mathbf{\beta}}%
))\,\{\dot{g}(x_{j},{\mathbf{\beta}})-\dot{g}(x_{i},{\mathbf{\beta}})\}.
\]
Assuming A0, by Lemma \ref{victor4.2}, we get
\begin{equation}
\mathbf{C}_{n}\longrightarrow\mathbf{C}=\mathrm{E}_{G_{0}}a_{i}\,I_{{T^{L}%
,F_{0}}}^{\prime}(y_{i}-g(x_{i},\mathbf{\beta}_{0})+g(x_{j},\mathbf{\beta}%
_{0}))\,\{\dot{g}(x_{j},\mathbf{\beta}_{0})-\dot{g}(x_{i},\mathbf{\beta}%
_{0})\},\text{a.s.} \label{TE5}%
\end{equation}
Using (\ref{TE1})-(\ref{TE5}), \ the expansion (\ref{TRexpansion}) guaranteed
by part (ii) of Theorem \ref{asslreg}, and the fact that \ by the U-statistics
projection Theorem $\{D_{n}\}$ converges to a normal distribution, we conclude
\ {that $\{\sqrt{n}\Vert$\textrm{$E$}$_{\widehat{F}_{n}}I_{T^{L},F_{0}%
}(y\mathbf{)}\Vert\}$ is bounded in probability.} Therefore, from (\ref{expb})
we get
\[
\sqrt{n}\{\mathbf{T}^{L}(\widehat{F}_{n})-\mathbf{T}^{L}(F_{0})\}=\sqrt
{n}\,\mathrm{E}_{\widehat{F}_{n}}I_{\mathbf{T}^{L},F_{0}}(y\mathbf{)+}%
o_{P}(1).
\]
This \ implies
\[
\sqrt{n}\{\mathbf{T}_{MM}^{L}(\widehat{F}_{n})-\mathbf{T}^{L}(F_{0}%
)\}=\sqrt{n}\,\mathrm{E}_{\widehat{F}_{n}}I_{\mathbf{T}_{MM}^{L},F_{0}%
}(y\mathbf{)+}o_{P}(1),
\]
and therefore \ (\ref{TLexpansion}) is satisfied with $I_{L}=I_{\mathbf{T}%
_{MM}^{L},F_{0}}.$ Finally \ (\ref{ifloc1}) follows \ from formula (44) of
\ Fasano et al. \cite{Fasano et al}. Part (ii) follows immediately from $e_{01}%
^{L}=0.\square$\bigskip

To prove Theorem \ref{locationbound}, the following result is required.

\begin{lemma}
\label{BDP1} Given $M$ and $\gamma>0,$ there exists $M^{\ast}$ such that
$\mathrm{P}_{F}(|y|\leq M)\geq1-\delta+\gamma$ implies $S^{L}(F)\leq M^{\ast
}.$
\end{lemma}

\begin{proof} It is enough to show that there exists $M^{\ast}$ such that
$S^{\ast }_{L}(F,0)\leq M^{\ast}$, where $S^{\ast}_{L}(F,\mu)$ is
the location version of the object defined by (\ref{Sestrella}) for
the regression case.

Let $M^{\ast}$ be such that $\rho_{0}^{L}(M/M^{\ast})<\gamma/2$. Suppose that
$S_{L}^{\ast}(F,0)>M^{\ast}$. By definition of $S_{L}^{\ast}(F,0)$,
\begin{equation}
\delta=\mathrm{E}_{F}\rho_{0}^{L}\left(  {y}/{S_{L}^{\ast}(F,0)}\right)
\ \label{SBP1}%
\end{equation}

On the other hand, let $A=\{|y|\leq M\}$. By hypothesis, \textrm{$P$}%
$_{F}(A)\geq1-\delta+\gamma$, and so
\[
\mathrm{E}_{F}\rho_{0}^{L}\left(  \frac{y}{S_{L}^{\ast}(F,0)}\right)
\leq\mathrm{E}_{F}\rho_{0}^{L}\left(  \frac{y}{M^{\ast}}\right)
\ \ \leq(\gamma/2)\mathrm{P}_{F}(A)+\mathrm{P}_{F}(A^{c})\leq\gamma
/2+\delta-\gamma\leq\delta-\gamma/2,
\]
contradicting (\ref{SBP1}). \end{proof}
\textsc{Proof of Theorem
\ref{locationbound}.} We will prove that, given $M$ and $\gamma>0,$
there exists $K$ such that $|T_{MM}^{L}(F)|\leq K$, for all $F$ with
$\mathrm{P}_{F}(|y|\leq M)\geq\min(1-\delta+\gamma,\delta+\gamma)$.
In fact , note that
\begin{equation}
\mathrm{E}_{F}\rho^{L}\left(  \frac{y-T_{MM}^{L}(F)}{S^{L}(F)}\right)
\leq\mathrm{E}_{F}\rho^{L}\left(  \frac{y-T_{S}^{L}(F)}{S^{L}(F)}\right)
\leq\mathrm{E}_{F}\rho_{0}^{L}\left(  \frac{y-T_{S}^{L}(F)}{S^{L}(F)}\right)
=\delta. \label{SBDP2}%
\end{equation}

Let $\ M^{\ast}$ be as in Lemma 1 and let $\ a$ so that $\rho^{L}(a/M^{\ast
}))(\delta+\gamma)=\delta+\gamma/2$. Put $K=M+a$ and observe that $|y|\leq M$
and $|T_{MM}^{L}(F)|>K$ imply that $|y-T_{MM}^{L}(F)|>a$. Suppose that
$|T_{MM}^{L}(F)|>K$ . Then
\[
\mathrm{E}_{F}\rho^{L}\left(  \frac{y-T_{MM}^{L}(F)}{S^{L}(F)}\right)
\geq\mathrm{E}_{F}\rho^{L}\left(  \frac{y-T_{MM}^{L}(F)}{M^{\ast}}\right)
\geq\mathrm{P}_{F}(A)\rho^{L}(a/M^{\ast}\geq\rho^{L}(a/M^{\ast})(\delta
+\gamma)\geq\delta+\gamma/2,
\]
contradicting (\ref{SBDP2}). $\square$

\ \ \ \ \ \ \ \

\end{document}